\documentclass[fleqn,preprint,3p,a4paper]{elsarticle}

\usepackage{amssymb}
\usepackage{amsmath}
\usepackage{amsthm}
\usepackage{dcolumn}
\usepackage{endnotes}
\usepackage{tabularx}
\usepackage[matrix,arrow]{xy}
\usepackage{wasysym}

\newtheorem{theorem}{Theorem}[section]
\newtheorem{proposition}[theorem]{Proposition}
\newtheorem{lemma}[theorem]{Lemma}
\newtheorem{corollary}[theorem]{Corollary}

\theoremstyle{definition}
\newtheorem{definition}[theorem]{Definition}
\newtheorem{example}[theorem]{Example}
\newtheorem{remark}[theorem]{Remark}
\newtheorem{question}[theorem]{Question}

\newcommand{\ir}{{\mathsf{Irr}}}

\newcommand{\mn}{\mathbb N}

\newcommand{\cl}{{\rm cl}}
\newcommand{\ii}{{\rm int}}

\newcommand{\ua}{\mathord{\uparrow}}
\newcommand{\da}{\mathord{\downarrow}}

\newcommand{\mk}{\mathord{\mathsf{K}}}
\newcommand{\wdd}{\mathord{\mathsf{WD}}}
\newcommand{\md}{\mathord{\mathsf{D}}}

\newcommand{\kf}{\mathord{\mathsf{RD}}}

\journal{Topology and its applications}

\begin{document}

\begin{frontmatter}



\title{First countability, $\omega$-well-filtered spaces and reflections\tnoteref{t1}}
\tnotetext[t1]{This research was supported by the National Natural Science Foundation of China (Nos. 11661057, 11361028,
61300153, 11671008, 11701500, 11626207); the Natural Science Foundation of Jiangxi Province , China (No. 20192ACBL20045); NSF Project of Jiangsu Province,
China (BK20170483); and NIE ACRF (RI 3/16 ZDS), Singapore}

\author[X. Xu]{Xiaoquan Xu\corref{mycorrespondingauthor}}
\cortext[mycorrespondingauthor]{Corresponding author}
\ead{xiqxu2002@163.com}
\address[X. Xu]{School of Mathematics and Statistics,
Minnan Normal University, Zhangzhou 363000, China}
\author[C. Shen]{Chong Shen}
\address[C. Shen]{School of Mathematics and Statistics,
Beijing Institute of Technology, Beijing 100081, China}
\ead{shenchong0520@163.com}
\author[X. Xi]{Xiaoyong Xi}
\ead{littlebrook@jsnu.edu.cn}
\address[X. Xi]{School of mathematics and Statistics,
Jiangsu Normal University, Xuzhou 221116, China}
\author[D. Zhao]{Dongsheng Zhao}
\address[D. Zhao]{Mathematics and Mathematics Education,
National Institute of Education Singapore, \\
Nanyang Technological University,
1 Nanyang Walk, Singapore 637616}
\ead{dongsheng.zhao@nie.edu.sg}

\begin{abstract}
  We first introduce and study two new classes of subsets in $T_0$ spaces - $\omega$-Rudin sets and $\omega$-well-filtered determined sets lying between the class of all closures of countable directed subsets and that of irreducible closed subsets, and two new types of spaces - $\omega$-$d$ spaces  and $\omega$-well-filtered spaces. We prove that an $\omega$-well-filtered $T_0$ space is locally compact if{}f  it is core compact. One immediate corollary is that every core compact well-filtered space is sober, answering Jia-Jung problem with a new method. We also prove that all irreducible closed subsets in a first countable $\omega$-well-filtered $T_0$ space are directed. Therefore, a first countable $T_0$ space $X$ is sober if{}f $X$ is well-filtered if{}f $X$ is an $\omega$-well-filtered $d$-space. Using $\omega$-well-filtered determined sets, we present a direct construction of the $\omega$-well-filtered reflections of $T_0$ spaces, and show that products of $\omega$-well-filtered spaces are $\omega$-well-filtered.
\end{abstract}




\end{frontmatter}


\section{Introduction}

In domain theory and non-Hausdorff topology, the $d$-spaces, well-filtered spaces and sober spaces form three of the most important classes of spaces (see [2-26]). In this paper, based on the topological version of Rudin's Lemma by Heckmann and Keimel \cite{Hofmann-Lawson}, we introduce and study two new classes of subsets in $T_0$ spaces - $\omega$-Rudin sets and $\omega$-well-filtered determined sets ($\omega$-$\wdd$ sets for short) lying between the class of all closures of countable directed subsets and that of irreducible closed subsets. We also introduce and investigate two new types of spaces - $\omega$-$d$ spaces  and $\omega$-well-filtered spaces. It will be  proved that an $\omega$-well-filtered $T_0$ space is locally compact if{}f  it is core compact. One immediate corollary is that every core compact well-filtered space is sober, giving a positive answer to Jia-Jung problem \cite{jia-2018}, which has been first answered by Lawson and Xi \cite{Lawson-Xi} using a different method. We also prove that every  irreducible closed subset of a first countable $\omega$-well-filtered $T_0$ space is directed. Therefore, a first countable $T_0$ space $X$ is sober if{}f $X$ is well-filtered if{}f $X$ is an $\omega$-well-filtered $d$-space.

It is well-known that the category of all sober spaces and that of $d$-spaces are reflective in the category of all $T_0$ spaces (see [7, 12, 23]).
Recently, following Ershov's method of constructing the $d$-completion of $T_0$ spaces, Shen, Xi, Xu and Zhao \cite {Shenchon} presented a construction of the well-filtered reflection of $T_0$ spaces.  In the current paper, using $\omega$-$\wdd$ sets, we present a direct construction of the $\omega$-well-filtered reflections of $T_0$ spaces, and show that products of $\omega$-well-filtered spaces are $\omega$-well-filtered. Some major properties of $\omega$-well-filtered reflections of $T_0$ spaces are also investigated.

\section{Preliminary}

In this section, we briefly recall some basic concepts and notations to be used in this  paper. Some known properties of irreducible sets and compact saturated sets are presented.

For a poset $P$ and $A\subseteq P$, let
$\mathord{\downarrow}A=\{x\in P: x\leq  a \mbox{ for some }
a\in A\}$ and $\mathord{\uparrow}A=\{x\in P: x\geq  a \mbox{
	for some } a\in A\}$. For  $x\in P$, we write
$\mathord{\downarrow}x$ for $\mathord{\downarrow}\{x\}$ and
$\mathord{\uparrow}x$ for $\mathord{\uparrow}\{x\}$. A subset $A$
is called a \emph{lower set} (resp., an \emph{upper set}) if
$A=\mathord{\downarrow}A$ (resp., $A=\mathord{\uparrow}A$). A nonempty subset $D$ of $P$ is \emph{directed} if every two
elements in $D$ have an upper bound in $D$. The set of all directed sets of $P$ is denoted by $\mathcal D(P)$. $P$ is called a
\emph{directed complete poset}, or \emph{dcpo} for short, if for any
$D\in \mathcal D(P)$, $\bigvee D$ exists in $P$. The set of all natural numbers with the usual ordering is denoted by $\mn$. Let $\omega$ denote the ordinal (also the cardinal number) of $\mn$ and $\omega_1$ the first uncountable ordinal.

The upper sets of a poset $Q$ form the (\emph{upper}) \emph{Alexandroff topology} $\alpha (Q)$ on $Q$. As in \cite{redbook}, the \emph{lower topology} on $Q$, generated
by the complements of the principal filters of $Q$, is denoted by $\omega (Q)$. A subset $U$ of $Q$ is \emph{Scott open} if
(i) $U=\mathord{\uparrow}U$ and (ii) for any directed subset $D$ for
which $\bigvee D$ exists, $\bigvee D\in U$ implies $D\cap
U\neq\emptyset$. All Scott open subsets of $Q$ form a topology,
and we call this topology  the \emph{Scott topology} on $Q$ and
denote it by $\sigma(P)$. The space $\Sigma ~\!Q=(Q, \sigma (Q))$ is called the \emph{Scott space} of $Q$. The common refinement $\sigma(Q)$ and $\omega (Q)$ is called the \emph{Lawson topology}
and is denoted by $\lambda (Q)$. The space $\Lambda (Q)=(Q, \lambda (Q))$ is called the \emph{Lawson space} of $Q$.

The category of all $T_0$ spaces is denoted by $\mathbf{Top}_0$. For any $X\in \mathbf{Top}_0$, $\leq_X$ denotes the \emph{specialization order} on $X$: $x\leq_X y$ if{}f $x\in \overline{\{y\}}$). In the following, when a $T_0$ space $X$ is considered as a poset, the order always refers to the specialization order if no other explanation. Let $\mathcal O(X)$ (resp., $\mathcal C(X)$) be the set of all open subsets (resp., closed subsets) of $X$. Define $\mathcal S_c(X)=\{\overline{{\{x\}}} : x\in X\}$ and $\mathcal D_c(X)=\{\overline{D} : D\in \mathcal D(X)\}$. A space $X$ is called a $d$-\emph{space} (or \emph{monotone convergence space}) if $X$ (with the specialization order) is a dcpo
 and $\mathcal O(X) \subseteq \sigma(X)$ (cf. \cite{redbook, Wyler}).

As in \cite{E_20182}, a space $X$ is \emph{locally hypercompact} if for each $x\in X$ and each open neighborhood $U$ of $x$, there is  $\ua F\in \mathbf{Fin}~X$ such that $x\in\ii\,\ua F\subseteq\ua F\subseteq U$. A space $X$ is called a $C$-\emph{space} if for each $x\in X$ and each open neighborhood $U$ of $x$, there is $u\in X$ such that $x\in\ii\,\ua u\subseteq\ua u\subseteq U$). A set $K\subseteq X$ is called \emph{supercompact} if for
any family $\{U_i : i\in I\}\subseteq \mathcal O(X)$, $K\subseteq \bigcup_{i\in I} U_i$  implies $K\subseteq U$ for some $i\in I$. It is easy to check that the supercompact saturated sets of $X$ are exactly the sets $\ua x$ with $x \in X$ (see \cite[Fact 2.2]{Klause-Heckmann}). It is well-known that $X$ is a $C$-space if{}f $\mathcal O(X)$ is a \emph{completely distributive} lattice (cf. \cite{E_2009}). A space $X$ is called \emph{core compact} if $\mathcal O(X)$ is a \emph{continuous lattice} (cf. \cite{redbook}).

For a $T_0$ space $X$ and a nonempty subset $A$ of $X$, $A$ is \emph{irreducible} if for any $\{F_1, F_2\}\subseteq \mathcal C(X)$, $A \subseteq F_1\cup F_2$ implies $A \subseteq F_1$ or $A \subseteq  F_2$.  Denote by $\ir(X)$ (resp., $\ir_c(X)$) the set of all irreducible (resp., irreducible closed) subsets of $X$. Clearly, every subset of $X$ that is directed under $\leq_X$ is irreducible. $X$ is called \emph{sober}, if for any  $F\in\ir_c(X)$, there is a unique point $a\in X$ such that $F=\overline{\{a\}}$.

The following two lemmas on irreducible sets are well-known.

\begin{lemma}\label{irrsubspace}
Let $X$ be a space and $Y$ a subspace of $X$. Then the following conditions are equivalent for a
subset $A\subseteq Y$:
\begin{enumerate}[\rm (1)]
	\item $A$ is an irreducible subset of $Y$.
	\item $A$ is an irreducible subset of $X$.
	\item ${\rm cl}_X A$ is an irreducible closed subset of $X$.
\end{enumerate}
\end{lemma}

\begin{lemma}\label{irrimage}
	If $f : X \longrightarrow Y$ is continuous and $A\in\ir (X)$, then $f(A)\in \ir (Y)$.
\end{lemma}

\begin{remark}\label{subspaceirr}  If $Y$ is a subspace of a space $X$ and $A\subseteq Y$, then by Lemma \ref{irrsubspace}, $\ir (Y)=\{B\in \ir(X) : B\subseteq Y\}\subseteq \ir (X)$ and  $\ir_c (Y)=\{B\in \ir(X) : B\in \mathcal C(Y)\}\subseteq \ir (X)$. If $Y\in \mathcal C(X)$, then $\ir_c(Y)=\{C\in \ir_c(X) : C\subseteq Y\}\subseteq \ir_c (X)$.
\end{remark}

\begin{lemma}\label{irrprod}\emph{(\cite{Shenchon})}
	Let	$X=\prod_{i\in I}X_i$ be the product space of $T_0$ spaces  $X_i (i\in I)$. If  $A$ is an irreducible subset of $X$, then $\cl_X(A)=\prod_{i\in I}\cl_{X_i}(p_i(A))$, where $p_i : X \longrightarrow X_i$ is the $i$th projection.
\end{lemma}

\begin{lemma}\label{prodirr}
	Let	$X=\prod_{i\in I}X_i$ be the product space of $T_0$ spaces  $X_i (i\in I)$  and $A_i\subseteq X_i (i\in I$). Then the following two conditions are equivalent:
\begin{enumerate}[\rm (1)]
	\item $\prod_{i\in I}A_i\in \ir (X)$.
	\item $A_i\in \ir (X_i)$ ($i\in I$).
\end{enumerate}
\end{lemma}
\begin{proof}  (1) $\Rightarrow$ (2): By Lemma \ref{irrimage}.

(2) $\Rightarrow$ (1): Let $A=\prod_{i\in I}A_i$. For $U, V\in \mathcal O(X)$, if $A\cap U\neq\emptyset\neq A\cap V$, then there exist $I_1, I_2\in I^{(<\omega)}$ and $(U_i, V_j)\in \mathcal O(X_i)\times \mathcal O(X_j)$ for all $(i, j)\in I_1\times I_2$ such that $\bigcap_{i\in I_1}p_i^{-1}(U_i)\subseteq U$, $\bigcap_{j\in I_2}p_j^{-1}(V_j)\subseteq V$ and $A\cap \bigcap_{i\in I_1}p_i^{-1}(U_i)\neq\emptyset\neq A\cap \bigcap_{j\in I_2}p_i^{-1}(V_j)$. Let $I_3=I_1\cup I_2$. Then $I_3$ is finite. For $i\in I_3\setminus I_1$ and $j\in I_3\setminus I_2$, let $U_i=X_i$ and $V_j=X_j$. Then for each $i\in I_3$, we have $A_i\cap U_i\neq\emptyset\neq A_i\cap V_i$, and whence $A_i\cap U_i\cap V_i\neq\emptyset$ by $A_i\in \ir (X_i)$. It follows that $A\cap \bigcap_{i\in I_1}p_i^{-1}(U_i)\cap \bigcap_{j\in I_2}p_i^{-1}(V_j)\neq\emptyset$, and consequently, $A\cap U\cap V\neq \emptyset$. Thus $A\in \ir (X)$.

\end{proof}

By Lemma \ref{irrprod} and Lemma \ref{prodirr}, we obtain the following corollary.

\begin{corollary}\label{irrcprod} Let $X=\prod_{i\in I}X_i$ be the product space of $T_0$ spaces  $X_i (i\in I)$. If  $A\in \ir_c(X)$, then $A=\prod_{i\in I}p_i(A)$ and $p_i(A)\in \ir_c (X_i)$ for each $i\in I$.
\end{corollary}

For any topological space $X$, $\mathcal G\subseteq 2^{X}$ and $A\subseteq X$, let $\Diamond_{\mathcal G} A=\{G\in \mathcal G : G\bigcap A\neq\emptyset\}$ and $\Box_{\mathcal G} A=\{G\in \mathcal G : G\subseteq  A\}$. The symbols $\Diamond_{\mathcal G} A$ and $\Box_{\mathcal G} A$ will be simply written as $\Diamond A$  and $\Box A$ respectively if there is no confusion. The \emph{lower Vietoris topology} on $\mathcal{G}$ is the topology that has $\{\Diamond U : U\in \mathcal O(X)\}$ as a subbase, and the resulting space is denoted by $P_H(\mathcal{G})$. If $\mathcal{G}\subseteq \ir (X)$, then $\{\Diamond_{\mathcal{G}} U : U\in \mathcal O(X)\}$ is a topology on $\mathcal{G}$. The space $P_H(\mathcal{C}(X)\setminus \{\emptyset\})$ is called the \emph{Hoare power space} or \emph{lower space} of $X$ and is denoted by $P_H(X)$ for short (cf. \cite{Schalk}). The \emph{upper Vietoris topology} on $\mathcal{G}$ is the topology that has $\{\Box_{\mathcal{G}} U : U\in \mathcal O(X)\}$ as a base, and the resulting space is denoted by $P_S(\mathcal{G})$.

\begin{remark} \label{eta continuous} Let $X$ be a $T_0$ space.
\begin{enumerate}[\rm (1)]
	\item If $\mathcal{S}_c(X)\subseteq \mathcal{G}$, then the specialization order on $P_H(\mathcal{G})$ is the order of set inclusion, and the \emph{canonical mapping} $\eta_{X}: X\longrightarrow P_H(\mathcal{G})$, given by $\eta_X(x)=\overline {\{x\}}$, is an order and topological embedding (cf. \cite{redbook, Jean-2013, Schalk}).
    \item The space $X^s=P_H(\ir_c(X))$ with the canonical mapping $\eta_{X}: X\longrightarrow X^s$ is the \emph{sobrification} of $X$ (cf. \cite{redbook, Jean-2013}).
\end{enumerate}
\end{remark}

A subset $A$ of a space $X$ is called \emph{saturated} if $A$ equals the intersection of all open sets containing it (equivalently, $A$ is an upper set in the specialization order). We shall use $\mathord{\mathsf{K}}(X)$ to
denote the set of all nonempty compact saturated subsets of $X$ and endow it with the \emph{Smyth preorder}, that is, for $K_1,K_2\in \mathord{\mathsf{K}}(X)$, $K_1\sqsubseteq K_2$ if{}f $K_2\subseteq K_1$. Let $\mathcal S^u(X)=\{\ua x : x\in X\}$. $X$ is called \emph{well-filtered} if it is $T_0$, and for any open set $U$ and filtered family $\mathcal{K}\subseteq \mathord{\mathsf{K}}(X)$, $\bigcap\mathcal{K}{\subseteq} U$ implies $K{\subseteq} U$ for some $K{\in}\mathcal{K}$.

For the well-filteredness of Scott topologies on dcpos, Xi and Lawson \cite{Xi-Lawson-2017}) gave the following interesting results.

\begin{lemma}\label{xi-lawson1}\emph{(\cite{Xi-Lawson-2017})} Let $P$ be a dcpo. If $(P, \lambda (P))$ is compact, then $(P, \sigma (P))$ is well-filtered.
\end{lemma}

\begin{corollary}\label{xi-lawson2}\emph{(\cite{Xi-Lawson-2017})} For any complete lattice $L$, $(L, \sigma (L))$ is well-filtered.
\end{corollary}

For any $T_0$ space $X$, the space $P_S(\mathord{\mathsf{K}}(X))$, denoted shortly by $P_S(X)$, is called the \emph{Smyth power space} or \emph{upper space} of $X$ (cf. \cite{Heckmann, Schalk}). It is easy to see that the specialization order on $P_S(X)$ is the Smyth order (that is, $\leq_{P_S(X)}=\sqsubseteq$). The \emph{canonical mapping} $\xi_X: X\longrightarrow P_S(X)$, $x\mapsto\ua x$, is an order and topological embedding (cf. \cite{Heckmann, Klause-Heckmann, Schalk}). Clearly, $P_S(\mathcal S^u(X))$ is a subspace of $P_S(X)$ and $X$ is homeomorphic to $P_S(\mathcal S^u(X))$.

 \begin{lemma}\label{X-Smyth-irr} Let $X$ be a $T_0$ space and $A\subseteq X$. Then the following conditions are equivalent:
 \begin{enumerate}[\rm (1)]
	\item $A\in\ir (X)$.
	\item $\xi_X(A)\in \ir (P_S(X))$.
\end{enumerate}
\end{lemma}
\begin{proof} (1) $\Rightarrow$ (2): By Lemma \ref{irrimage}.

(2) $\Rightarrow$ (1): Suppose that $B, C\in \mathcal C(X)$ such that $A\subseteq B\cup C$. Then $\Diamond B, \Diamond C\in \mathcal C(P_S(X))$ and $\xi_X(A)\subseteq \Diamond B\cup \Diamond C$, and hence $\xi_X(A)\subseteq \Diamond B$ or $\xi_X(A)\subseteq \Diamond C$ by $\xi_X(A)\in \ir (P_S(X))$. It follows that $A\subseteq B$ or $A\subseteq C$. Thus $A\in\ir (X)$.

\end{proof}

\begin{remark}\label{meet-in-Smyth} Let $X$ be a $T_0$ space and $\mathcal A\subseteq \mk (X)$. Then $\bigcap \mathcal A=\bigcap \overline{\mathcal A}$, here the closure of $\mathcal A$ is taken in $P_S(X)$. Clearly, $\bigcap \overline{\mathcal A}\subseteq\bigcap \mathcal A$. On the other hand, for $K\in \overline{\mathcal A}$ and $U\in \mathcal O(X)$ with $K\subseteq U$ (that is, $K\in \Box U$), we have $\mathcal A\bigcap\Box U\neq\emptyset$, and hence there is a $K_U\in \mathcal A\bigcap\Box U$. Therefore, $K=\bigcap \{U\in \mathcal O(X) : K\subseteq U\}\supseteq\bigcap \{K_U : U\in \mathcal O(X) \mbox{ and } K\subseteq U\}\supseteq\bigcap \mathcal A$. It follows that $\bigcap \overline{\mathcal A}\supseteq\bigcap \mathcal A$. Thus $\bigcap \mathcal A=\bigcap \overline{\mathcal A}$.
\end{remark}

\begin{lemma}\label{Kmeet}\emph{(\cite{redbook})} For a nonempty family $\{K_i : i\in I\}\subseteq \mk (X)$, $\bigvee_{i\in I} K_i$ exists in $\mk (X)$ if{}f~$\bigcap_{i\in I} K_i\in \mk (X)$. In this case $\bigvee_{i\in I} K_i=\bigcap_{i\in I} K_i$.
\end{lemma}

\begin{lemma}\label{K union} \emph{(\cite{jia-Jung-2016, Schalk})}  Let $X$ be a $T_0$ space. If $\mathcal K\in\mk(P_S(X))$, then  $\bigcup \mathcal K\in\mk(X)$.
\end{lemma}

\begin{corollary}\label{Smythunioncont} \emph{(\cite{jia-Jung-2016, Schalk})}  For any $T_0$ space $X$, the mapping $\bigcup : P_S(P_S(X)) \longrightarrow P_S(X)$, $\mathcal K\mapsto \bigcup \mathcal K$, is continuous.
\end{corollary}
\begin{proof} For $\mathcal K\in\mk(P_S(X))$, $\bigcup \mathcal K=\bigcup \mathcal K\in\mk(X)$ by Lemma \ref{K union}. For $U\in \mathcal O(X)$, we have $\bigcup^{-1}(\Box U)=\{\mathcal K\in \mk (P_S(X)) : \bigcup \mathcal K\in \Box U\}=\{\mathcal K\in \mk (P_S(X)) : \mathcal K\subseteq \Box U\}=\eta_{P_S(X)}^{-1}(\Box (\Box U))\in \mathcal O(P_S(P_S(X)))$. Thus  $\bigcup : P_S(P_S(X)) \longrightarrow P_S(X)$ is continuous.
\end{proof}

\section{$\omega$-$d$-spaces and $\omega$-well-filtered spaces}

For a $T_0$ space $X$, let $\mathcal D^\omega(X)=\{D\subseteq X : D \mbox{ is countable and directed}\}$ and $\mathcal D^\omega_c(X)=\{\overline{D} : D\in \mathcal D^\omega(X)\}$.

\begin{definition}
	A poset $P$ is called an \emph{$\omega$-dcpo}, if for any $D\in \mathcal D^\omega(X)$, $\bigvee D$ exists.
\end{definition}

\begin{lemma}\label{cchain}
Let $P$ be a poset and $D$ a countable directed subset of $P$. Then there exists a countable chain $C\subseteq D$ such that $D=\da C$. Hence, $\bigvee C$ exists and $\bigvee C=\bigvee D$ whenever $\bigvee D$ exists.
\end{lemma}
\begin{proof}
	If $|D|<\omega$, then $D$ contains a largest element $d$, so let $C=\{d\}$, which satisfies the requirement.
	
	Now assume $|D|=\omega$ and let $D=\{d_n:n<\omega\}$. We use induction on $n\in\omega$ to  define $C=\{c_n:n<\omega\}$.
	More precisely, let $c_0=d_0$ and let $c_{n+1}$ be an upper bound of $\{d_{n+1},c_0, c_1,c_2\ldots,c_n\}$ in $D$. It is clear that $C$ is a chain and $D=\da C$.
\end{proof}

\begin{corollary}
	A poset $P$ is an $\omega$-dcpo if{}f for any countable chain $C$ of $P$, $\bigvee C$ exists.
\end{corollary}

\begin{definition}
	Let $P$ be a poset. A subset $U$ of $P$ is called \emph{$\omega$-Scott open} if (i) $U=\ua U$, and (ii) for any countable directed set $D$, $\bigvee D\in U$ implies that $D\cap U\neq\emptyset$. All $\omega$-Scott open sets form a topology on $P$,  denoted by $\sigma_\omega(P)$ and called the \emph{$\omega$-Scott topology}. The space $\Sigma_{\omega}P=(P,\sigma_{\omega}(P))$ is called the \emph{$\omega$-Scott space} of $P$.
\end{definition}

\begin{remark}
(1)	By Lemma \ref{cchain}, $U\in\sigma_{\omega}(P)$ if{}f $U=\ua U$ and for any countable chain $C$, $\bigvee C\in U$ implies $C\cap U\neq\emptyset$.

(2) Clearly, $\sigma(P)\subseteq\sigma_{\omega}(P)$. The converse need not be true, see Example \ref{omegat}.
\end{remark}

$\mathcal D^\omega(X)=\{D\subseteq X : D \mbox{ is countable directed}\}$ and $\mathcal D^\omega_c(X)=\{\overline{D} : D\in \mathcal D^\omega(X)\}$

\begin{definition}
	A $T_0$ space $X$ is called an $\omega$-$d$-space (or an \emph{$\omega$-monotone convergence space}) if for any $D\in\mathcal D^\omega(X)$, the closure of $D$ has a generic point, equivalently, if $\mathcal D_c^\omega(X)=\mathcal S_c(X)$.
\end{definition}

\begin{proposition}\label{d-spacecharac1} For a $T_0$ space $X$, the following conditions are equivalent:
\begin{enumerate}[\rm (1) ]
	        \item $X$ is an $\omega$-$d$-space.
            \item $X$ (with the specialization order $\leq_X$) is an $\omega$-dcpo and $\mathcal O(X)\subseteq \sigma_{\omega}(X)$.
            \item  For any $D\in \mathcal D^\omega(X)$ and $U\in \mathcal O(X)$, $\bigcap\limits\limits_{d\in D}\ua d\subseteq U$ implies $\ua d \subseteq U$ \emph{(}i.e., $d\in U$\emph{)} for some $d\in D$.
            \item  For any $D\in \mathcal D^\omega(X)$ and $A\in \mathcal C(X)$, if $D\subseteq A$, then $A\cap\bigcap\limits_{d\in D}\ua d\neq\emptyset$.
            \item  For any $D\in \mathcal D^\omega(X)$ and $A\in \ir_c(X)$, if $D\subseteq A$, then $A\cap\bigcap\limits_{d\in D}\ua d\neq\emptyset$.
            \item  For any $D\in \mathcal D^\omega(X)$, $\overline{D}\cap\bigcap\limits_{d\in D}\ua d\neq\emptyset$.
\end{enumerate}
\end{proposition}
\begin{proof} (1) $\Leftrightarrow$ (2): Clearly, (2) $\Rightarrow$ (1). Conversely, if condition (1) holds, then for each $D\in \mathcal D^\omega(X)$ and $A\in \mathcal C(X)$ with $D\subseteq A$, there is $x\in X$ such that $\overline{D}=\overline{\{x\}}$, and consequently, $\bigvee D=x$ and $\bigvee D \in A$ since $\overline{D}\subseteq A$. Thus $X$ is an $\omega$-dcpo and $\mathcal O(X)\subseteq \sigma_\omega(X)$.

(2) $\Rightarrow$ (3): By condition (2), $\ua \bigvee D=\bigcap\limits_{d\in D}\ua d\subseteq U\in \sigma_\omega (X)$. Therefore, $\bigvee D\in U$, and whence
$d\in U$ for some $d\in D$.

(3) $\Rightarrow$ (4): If $A\cap\bigcap\limits_{d\in D}\ua d=\emptyset$, then $\bigcap\limits_{d\in D}\ua d \subseteq X\setminus A$. By condition (3), $\ua d\subseteq X\setminus A$ for some $d\in D$, which is in contradiction with $D\subseteq A$.

(4) $\Rightarrow$ (5) $\Rightarrow$ (6): Trivial.

(6) $\Rightarrow$ (1): For each $D\in \mathcal D^\omega(X)$ and $A\in \mathcal C(X)$ with $D\subseteq A$, by condition (6), $\overline{D}\cap\bigcap\limits_{d\in D}\ua d\neq\emptyset$. Select an $x\in \overline{D}\cap\bigcap\limits_{d\in D}\ua d$. Then $D\subseteq \da x\subseteq \overline{D}$, and hence $\overline{D}=\da x$. Thus $X$ is an $\omega$-$d$-space.
\end{proof}

$P$ is said to be \emph{$\omega$-Noetherian} if it satisfies the $\omega$-\emph{ascending chain condition} ($\omega$-$\mathrm{ACC}$ for short): every countable ascending chain has a greatest member. By Lemma \ref{cchain}, $P$ is $\omega$-Noetherian if{}f every countable directed set of $P$ has a largest element.

For a poset $P$ and $x\in P$, $x$ is called an \emph{$\omega$-compact} element of $P$, written as $x\ll_\omega x$, if for each countable directed subset $D$ of $P$ with $\bigvee D$ exists, $x\leq \bigvee D$ implies $x\leq d$ for some $d\in D$.

\begin{proposition}\label{alphaOmegad}
	For a poset $P$, the following conditions are equivalent:
	\begin{enumerate}[\rm (1)]
		\item $(P,\alpha(P))$ is an $\omega$-$d$-space.
        \item $P$ is an $\omega$-Noetherian poset.
		\item $P$ is an $\omega$-dcpo such that every element of $P$ is $\omega$-compact \emph{(}i.e., $x\ll_\omega x$ for all $x\in P$ \emph{)}.
		\item $P$ is an $\omega$-dcpo such that $\alpha(P)=\sigma_\omega(P)$.
	\end{enumerate}
\end{proposition}

\begin{proof} (1) $\Rightarrow$ (2): Suppose that $D$ is a countable directed set of $P$. Then $\da D$ is a closed subset in $(P,\alpha(P))$. Since $(P,\alpha(P))$ is an $\omega$-$d$-space, by Proposition \ref{d-spacecharac1}, $\bigcap_{d\in D}\ua d\cap \da D\neq \emptyset$. Clearly, $x\in \bigcap_{d\in D}\ua d\cap \da D$ if{}f $x$ is the greatest element of $P$.

(2) $\Rightarrow$ (3)and (4) $\Rightarrow$ (1): Trivial.

(3) $\Rightarrow$ (4): For any $x\in X$, since $x\ll_\omega x$, we have $\ua x\in \sigma_\omega(P)$. Therefore, $\alpha(P)=\sigma_\omega(P)$.

\end{proof}

\begin{definition}
	A $T_0$ space $X$ is called \emph{$\omega$-well-filtered}, if for any countable filtered family $\{K_i:i<\omega\}\subseteq \mk (X)$ and  $U\in\mathcal O(X)$, it satisfies
		$$\bigcap_{i<\omega}K_i\subseteq U \ \Rightarrow \  \exists i_0<\omega, K_{i_0}\subseteq U.$$
\end{definition}

By Lemma \ref{cchain}, we have the following result.

\begin{proposition}
	 A $T_0$ space $X$ is $\omega$-well-filtered if{}f for any countable descending chain $K_0\supseteq K_1\supseteq K_2\supseteq\ldots\supseteq K_n\supseteq\ldots$ of compact saturated subsets of $X$ and $U\in\mathcal O(X)$, the following implication holds:
	$$\bigcap_{i<\omega}K_i\subseteq U\ \Rightarrow\ \exists i_0<\omega, \ K_{i_0}\subseteq U.$$
\end{proposition}

By Proposition \ref{d-spacecharac1}, we get the following result.

\begin{proposition}\label{wfismc}
	Every $\omega$-well-filtered space is an $\omega$-$d$-space.
\end{proposition}

The following result is well-known.

\begin{theorem}\label{SoberLC=CoreC}\emph{(\cite{redbook, Jean-2013, Kou})}  For a $T_0$ space $X$, the following conditions are equivalent:
\begin{enumerate}[\rm (1)]
	\item $X$ locally compact and sober.
	\item $X$ is locally compact and well-filtered.
	\item $X$ is core compact and sober.
\end{enumerate}
\end{theorem}

The above theorem will be strengthened (see Theorem \ref{SoberLC=CoreCNew}). First, we give the following result.

\begin{theorem}\label{omegaWFcorecompt-locCompt} For an $\omega$-well-filtered $T_0$ space $X$, $X$ is locally compact if{}f $X$ is core compact.
\end{theorem}

\begin{proof}  Suppose that $X$ is core compact. For $x\in X$ and $U\in \mathcal O(X)$ with $x\in U$, since $X$ is core compact, there is a sequence $\{U_\infty, ..., U_n, ..., U_2, U_1, U_0\}\subseteq \mathcal O(X)$ such that
$$x\in U_\infty\ll ...\ll U_n\ll ...\ll U_2\ll U_1\ll U_0=U.$$
We show that $K=\bigcap\limits_{i\in \mathbb{N}}U_i$ is compact. Suppose that $\{W_i:i\in I\}\subseteq \mathcal O(X)$ is an open cover of $K$. Let $W=\bigcup_{i\in I}W_i$. If $U_n \not\subseteq W$ for all $n\in \mathbb{N}$, then by Lemma \ref{t Rudin}, there is a minimal irreducible closed subset $C\subseteq X\setminus W$ such that $U_n\cap C\neq\emptyset$. For each $n\in \mathbb{N}$, select an $x_n\in U_n\cap C$ and let $H_n=\{x_m : n\leq m\}$. Now we prove that $\ua H_n\in \mk (X)$ for all $n\in \mathbb{N}$. Suppose that $\{V_d : d\in D\}\subseteq \mathcal O(X)$ is a directed open cover of $\ua H_n$.

     {(c1)} If for some $d_1\in D$, $H_n\cap (X\setminus V_{d_1})=H_n\setminus V_{d_1}$ is finite, then $H_n\setminus V_{d_1}\subseteq V_{d_2}$ for some $d_2\in D$ because $H_n\subseteq \bigcup\limits_{d\in D}V_d$ and $\{V_d : d\in D\}\subseteq \mathcal O(X)$ is directed. By the directness of $\{V_d : d\in D\}$ again, $V_{d_1}\cup V_{d_2}\subseteq V_{d_3}$ for some $d_3\in D$. Then $H_n\subseteq V_{d_3}$.

      {(c2)} If for all $d\in D$, $H_n\cap (X\setminus V_{d})$ is infinite, then $U_n\cap C\cap (X\setminus V_d)\neq\emptyset$ since $U_1\supseteq U_2\supseteq...\supseteq U_n\supseteq ...$, and whence $C\cap (Y\setminus V_d)=C$ for all $d\in D$ by the minimality of $C$. Therefore, $H_n\subseteq C\subseteq \bigcap_{d\in D}(X\setminus V_d)=X\setminus \bigcup_{d\in D}V_d$, which is a contradiction with $H_n\subseteq \bigcup_{d\in D}V_d$.

By (c1) and (c2), $\ua H_n\in  \mk (X)$. Clearly, $\{\ua H_n : n\in N\}\subseteq \mk (X)$ is countable filtered, and whence $H=\bigcap_{n\in N}\ua H_n\in \mk (X)$ by the $\omega$-well-filteredness of $X$. It follows that $\emptyset \neq H\subseteq C\cap \bigcap_{n\in \mathbb{N}}U_n\subseteq (X\setminus W)\cap W=\emptyset$ (note that $U_1\supseteq U_2\supseteq ... \supseteq U_n\supseteq ...$), a contradiction.

Therefore, $U_{n_0}\subseteq W=\bigcup_{i\in I}W_i$ for some $n_0\in \mathbb{N}$. By $U_{n_0+1}\ll U_{n_0}$, $K\subseteq U_{n_0+1}\subseteq \bigcup_{i\in J}W_i$ for some $J\in I^{<\omega}$. It follows that $K\in \mk (X)$ and $x\in U_\infty \subseteq K\subseteq U$. Thus $X$ is locally compact.
\end{proof}

\begin{corollary}\label{WFcorcomp-sober}  A well-filtered  $T_0$ space is locally compact if{}f it is core compact.
\end{corollary}

By Theorem \ref{SoberLC=CoreC} and Theorem \ref{omegaWFcorecompt-locCompt}, we reobtain  the following result, which was first proved by Lawson and Xi \cite{Lawson-Xi} using a different method.

\begin{corollary}\label{WFcorcomp-sober}  Every well-filtered core compact $T_0$ space is sober.
\end{corollary}

Corollary \ref{WFcorcomp-sober} gives a positive answer to Jia-Jung problem \cite{jia-2018} (see \cite[Question 2.5.19]{jia-2018}). The Theorem \ref{SoberLC=CoreC} can be strengthened into the following one.

\begin{theorem}\label{SoberLC=CoreCNew}  For a $T_0$ space $X$, the following conditions are equivalent:
\begin{enumerate}[\rm (1)]
	\item $X$ locally compact and sober.
	\item $X$ is locally compact and well-filtered.
	\item $X$ is core-compact and sober.
    \item $X$ is core compact and well-filtered.
\end{enumerate}
\end{theorem}

Rudin's Lemma \cite{Rudin} is a very useful tool in domain theory and non-Hausdorff topology (see [2-10, 13, 22]). In \cite{Klause-Heckmann}, Heckman and Keimel presented the following topological variant of Rudin's Lemma.

\begin{lemma}\label{t Rudin} \emph{(Topological Rudin's Lemma)} Let $X$ be a topological space and $\mathcal{A}$ an
irreducible subset of the Smyth power space $P_S(X)$. Then every closed set $C {\subseteq} X$  that
meets all members of $\mathcal{A}$ contains an minimal irreducible closed subset $A$ that still meets all
members of $\mathcal{A}$.
\end{lemma}

In the following, using the topological Rudin's Lemma, we prove that a $T_0$ space $X$ is  $\omega$-well-filtered if{}f the Smyth power space of $X$ is $\omega$-well-filtered. The corresponding results for well-filteredness were given in \cite{xi-zhao-MSCS-well-filtered}\cite{xu-xi-zhao-19}.

\begin{theorem}\label{Smythomegawf}
	For a $T_0$ space, the following conditions are equivalent:
\begin{enumerate}[\rm (1)]
		\item $X$ is $\omega$-well-filtered.
        \item $P_S(X)$ is an ${\omega}$-$d$-space.
        \item $P_S(X)$ is $\omega$-well-filtered.
\end{enumerate}
\end{theorem}
\begin{proof}

(1) $\Rightarrow$ (2): Suppose that $X$ is an  $\omega$-well-filtered space. For any countable $\mathcal K\in \mk (X)$, by the $\omega$-well-filteredness of $X$,  $\bigcap \mathcal K \in \mk (X)$. Therefore,  by Lemma \ref{Kmeet}, $\mk (X)$ is an ${\omega}$-dcpo. Clearly, by the $\omega$-well-filteredness of $X$,  $\Box U\in \sigma_\omega (\mk (X))$ for any $U\in O(X)$. Thus $P_S (X)$ is an ${\omega}$-$d$-space.

(2) $\Rightarrow$ (3): Suppose that $\{\mathcal K_n : n<\omega\}\subseteq \mk(P_S(X))$ is countable filtered, $\mathcal U\in \mathcal O(P_S(X))$, and $\bigcap\limits_{n<\omega} \mathcal K_n \subseteq \mathcal U$. If $\mathcal K_n\not\subseteq \mathcal U$ for all $n<\omega$, then by Lemma \ref{t Rudin}, $\mk (X)\setminus \mathcal U$ contains an irreducible closed subset $\mathcal A$ that still meets all $\mathcal K_n$ ($n<\omega$). For each $n<\omega$, let $K_n=\bigcup \mathcal \ua_{\mk (X)} (\mathcal A\bigcap \mathcal K_n)$ ($=\bigcup (\mathcal A\bigcap \mathcal K_n$)). Then by Lemma \ref{K union}, $\{K_n : n<\omega\}\subseteq \mk (X)$ is countable filtered, and $K_n\in \mathcal A$ for all $n<\omega$ since $\mathcal A=\da_{\mk (X)}\mathcal A$. Let $K=\bigcap\limits_{n<\omega} K_n$. Then $K\in \mk (X)$ and $K=\bigvee_{\mk (X)} \{K_n : n<\omega\}\in \mathcal A$ by Lemma \ref{Kmeet} and condition (2). We claim that $K\in \bigcap\limits_{n<\omega}\mathcal K_n$. Suppose, on the contrary, that $K\not\in \bigcap\limits_{n<\omega}\mathcal K_n$. Then there is a $n_0<\omega$ such that $K\not\in \mathcal K_{n_0}$. Select a $G\in \mathcal A\bigcap \mathcal K_{n_0}$. Then $K\not\subseteq G$ (otherwise, $K\in \ua_{\mk (X)}\mathcal K_{n_0}=\mathcal K_{n_0}$, being a contradiction with $K\not\in \mathcal K_{n_0}$), and hence there is a $g\in K\setminus G$. It follows that $g\in K_n=\bigcup (\mathcal A\bigcap \mathcal K_n)$ for all $n<\omega$ and $G\not\in \Diamond_{\mk (K)}\overline{\{g\}}$. For each $n<\omega$, by $g\in K_n=\bigcup (\mathcal A\bigcap \mathcal K_n)$, there is a $K_n^g\in \mathcal A\bigcap \mathcal K_n$ such that $g\in K_n^g$, and hence $K_n^g\in\Diamond_{\mk (K)}\overline{\{g\}}\bigcap\mathcal A\bigcap \mathcal K_n$. Thus $\Diamond_{\mk (K)}\overline{\{g\}}\bigcap\mathcal A\bigcap \mathcal K_n\neq\emptyset$ for all $n<\omega$. By the minimality of $\mathcal A$, we have $\mathcal A=\Diamond_{\mk (K)}\overline{\{g\}}\bigcap\mathcal A$, and consequently, $G\in \mathcal A\bigcap \mathcal K_{n_0}=\Diamond_{\mk (K)}\overline{\{g\}}\bigcap\mathcal A\bigcap \mathcal K_{n_0}$, which is a contradiction with $G\not\in \Diamond_{\mk (K)}\overline{\{g\}}$. Thus $K\in \bigcap\limits_{n<\omega}\mathcal K_n\subseteq \mathcal U\subseteq \mk (X)\setminus \mathcal A$, being a contradiction with $K\in \mathcal A$. Therefore, $P_S(X)$  is $\omega$-well-filtered.

(3) $\Rightarrow$ (1): Suppose that $\mathcal K\subseteq \mathord{\mathsf K}(X)$ is countable filtered, $U\in \mathcal O(X)$, and $\bigcap \mathcal K \subseteq U$. Let $\widetilde{\mathcal K}=\{\ua_{\mk (X)}K : K\in \mathcal K\}$. Then $\widetilde{\mathcal K}\subseteq \mk (P_S(X))$ is countable filtered and $\bigcap \widetilde{\mathcal K} \subseteq \Box U$. By the $\omega$-well-filteredness of $P_S(X)$, there is a $K\in \mathcal K$ such that $\ua_{\mk (X)}K\subseteq \Box U$, and whence $K\subseteq U$, proving that $X$ is $\omega$-well-filtered.
\end{proof}

\section{First countable $\omega$-well-filtered spaces}

In this section, we show that in a first countable $\omega$-well-filtered $T_0$ space $X$, all irreducible closed subsets of $X$ are directed. Therefore, every first countable $\omega$-well-filtered $d$-space (in particular, every first countable well-filtered $T_0$ space) is sober.

\begin{theorem}\label{maintheorem1} Let $X$ be a first countable $\omega$-well-filtered $T_0$ space and $A\in \ir (X)$. Then $\overline{A}$ is directed.
\end{theorem}

\begin{proof}  For each $x\in X$, since $X$ is first countable, there is an open neighborhood base $\{U_n(x) : n\in \mn\}$ of $x$ such that

$$U_1{(x)}\supseteq U_2{(x)}\supseteq\ldots\supseteq U_k{(x)}\supseteq\ldots, $$

\noindent that is,  $\{U_n(x) : n\in \mn\}$ is a decreasing sequence of open subsets.

Suppose that $A\in \ir (X)$. We show that $\overline{A}$ is directed. Let $a_1,a_2\in \overline{A}$. It needs to show $\ua a_1\cap \ua a_2\cap \overline{A}\neq\emptyset$.
 Since $a_1,a_2\in \overline{A}$, $A\cap U_1(a_1)\neq\emptyset \neq A\cap U_1(a_2)$, and hence $A\cap U_1{(a_1)}\cap U_1{(a_2)}\neq\emptyset$ by the irreducibility of $A$. Choose $c_1\in  U_1{(a_1)}\cap U_1{(a_2)}\cap A$.
Now suppose we already have a set $\{c_1,c_2,\ldots, c_n\}$ such that
for each $2\leq i\leq n$,
$$c_{i}\in U_{i}(c_1)\cap U_{i}(c_2)\cap \ldots\cap U_{i}(c_i)\cap U_{i}(a_1)\cap U_{i}(a_2)\cap A.$$
Note that above condition implies that for any positive integer $k$,
$$U_{k}(c_1)\cap U_k(c_2)\cap \ldots\cap U_k(c_n)\cap U_k(a_1)\cap U_k(a_2)\cap A\neq \emptyset.$$
So we can choose $c_{n+1}\in U_{n+1}(c_1)\cap U_{n+1}(c_2)\cap \ldots\cap U_{n+1}(c_n)\cap U_{n+1}(a_1)\cap U_{n+1}(a_2)\cap A\neq\emptyset$.
By induction, we can obtain a set $\{c_n: n\in\mn\}$.

Let $K_n=\ua\{c_k:k\geq n\}$ for each $n\in\mn$.

{Claim 1:} $\forall n\in\mn$, $K_n$ is compact.

Suppose $\{V_i:i\in I\}$ is an open cover of $K_n$, i.e., $K_n\subseteq \bigcup_{i\in I}V_i$. Then there is $i_0\in I$ such that $c_n\in V_{i_0}$, and thus there is $m\geq n$ such that $c_n\in U_m(c_n)\subseteq V_{i_0}$. It follows that $c_k\in U_m(c_n)\subseteq V_{i_0}$ for all $k\geq m$. Thus $\{c_k: k\geq m\}\subseteq V_{i_0}$. For each $c_k$, where $n\leq k\leq m$, choose a $V_{i_k}$ such that $c_k\in V_{i_k}$.
Then the finite family $\{V_{i_k}:n\leq k\leq m\}\cup \{V_{i_0}\}$  covers $K_n$. So $K_n$ is compact.

{Claim 2:} $\bigcap_{n\in \mn}K_n\cap \overline{A}\neq \emptyset$.

Assume  $\bigcap_{n\in \mn}K_n\cap \overline{A}= \emptyset$. Then $\bigcap_{n\in \mn}K_n\subseteq X\setminus\overline{A}$.
Since $\{K_n: n\in\mn\}$ is a countable filtered family of compact saturated set in $X$ and $X$ is $\omega$-well-filtered, there exists $n_0\in\mn$ such that $K_{n_0}\subseteq X\setminus \overline{A}$, a contradiction.

{Claim 3:} $\ua a_1\cap \ua a_2\cap\overline{A}\neq\emptyset$.

Note that for each $n\in\mn$, $K_n=\ua \{c_k: k\geq n\}\subseteq U_{n}(a_1)\cap U_{n}(a_2)$.  This implies that
$$\emptyset\neq\bigcap_{n\in\mn}K_n\cap \overline{A}\subseteq\bigcap_{n\in\mn}U_n(a_1)\cap \bigcap_{n\in\mn}U_n(a_2)\cap \overline{A}=\ua a_1\cap \ua a_2\cap \overline{A}.$$
Thus $\ua a_1\cap \ua a_2\cap \overline{A}\neq\emptyset$, and whence there is $a_3\in \ua a_1\cap \ua a_2\cap \overline{A}$, that is, $a_3\in \overline{A}$ such that $a_1\leq a_3$ and $a_2\leq a_3$. Therefore, $\overline{A}$ is directed.
\end{proof}

\begin{theorem}\label{maintheorem2}
	Let $X$ be a first countable $T_0$ space. Then the following conditions are equivalent:
\begin{enumerate}[\rm (1)]
	\item $X$ is a sober space.
	\item $X$ is a well-filtered space.
	\item $X$ is an $\omega$-well-filtered $d$-space.
\end{enumerate}
\end{theorem}
\begin{proof}

(1) $\Rightarrow$ (2) $\Rightarrow$ (3): Trivial.

(3) $\Rightarrow$ (1): Let $A\in \ir_c(X)$. Then by Theorem \ref{maintheorem1}, $A\in \mathcal D(X)$. Since $X$ is a $d$-space, $\bigvee A\in A$, and hence $A=\overline{\{\bigvee A\}}$. Thus $X$ is sober.
\end{proof}

\begin{example}\label{omegat}
	Let $L$ be the complete chain $[0,\omega_1]$. Consider the space $\Sigma_{\omega}L=(L,\sigma_{\omega}(L))$.
\begin{enumerate}[\rm (1)]
	\item It is a first countable $\omega$-well-filtered space.
	\item It is not a $d$-space. In fact, we have that $\{\omega_1\}\in\sigma_{\omega}(L)$ but $\{\omega_1\}\notin\sigma(L)$.
  \end{enumerate}
  Therefore,  $\Sigma_{\omega}L$ is not well-filtered, and hence non-sober. So in Theorem \ref{maintheorem2}, condition (3) cannot be weakened to the condition that $X$ is only an $\omega$-well-filtered space.	
\end{example}

\begin{example} Let $L$ be the complete lattice constructed by Isbell \cite{isbell}. Then $\Sigma L$ is non-sober. By Corollary \ref{xi-lawson2} and Theorem \ref{maintheorem2}, $\Sigma L$ is well-filtered but not first countable.
\end{example}

Recently, Jung\footnote{A. Jung, Four dcpos, a theorem, and an open problem, an academic report at National Institute of Education, Singapore, 1 February, 2019} asked whether there is a countable complete lattice whose Scott space is non-sober. If there is such a countable complete lattice $L$, then $(L, \sigma (L))$ cannot be first countable (see Corollary \ref{Complete-case}).

By Lemma \ref{xi-lawson1}, Corollary \ref{xi-lawson2} and Theorem \ref{maintheorem2}, we get the following results.

\begin{corollary}\label{dcpo-case1} For a dcpo $P$, if $(P, \sigma (P))$ is first countable and $(P, \lambda (P))$ is compact, then $(P, \sigma (P))$ is sober.
\end{corollary}

\begin{corollary}\label{dcpo-case2} For a dcpo $P$, if $(P, \sigma (P))$ is a first countable $\omega$-well-filtered space, then it is sober.
\end{corollary}

\begin{corollary}\label{Complete-case} For a complete lattice $L$, if $(L, \sigma (L))$ is first countable, then it is sober.
\end{corollary}

The reader may wonder whether we can answer Jung's question positively by showing that the Scott topology on every countable complete lattice is first countable. Unfortunately, the following example crashes this hope.

\begin{example}
	Let $L=\{\bot\}\cup(\mn \times \mn)\cup\{\top\}$ and define a partial order $\leq$ on $L$ as follows:
	\begin{itemize}
		\item [(i)] $\forall (n,m)\in \mn\times\mn$, $\bot\leq (n,m) \leq\top$;
		\item [(ii)] $\forall (n_1,m_1), (n_2,m_2)\in\mn\times\mn$, $(n_1,m_1)\leq(n_2,m_2)$ iff $n_1=n_2$ and $m_1\leq m_2$.
	\end{itemize}

We show that $(L, \sigma (L))$ does not have any countable base at $\top$. Assume, on the contrary,  there exists a countable base $\{U_n : n\in\mn\}$ at $\top$. Then for each $n\in\mn$, as
$$\bigvee (\{n\}\times\mn)=\top\in U_n,$$
there exists $m_n\in\mn$ such that $(n,m_n)\in U_n$.
Let $U=\bigcup_{n\in\mn }\ua (n,m_n+1)$. Then $U\in \sigma (L)$. But for each $n\in N$, $(n,m_{n})\in U_n\setminus U$, which contradicts that  $\{U_n:n\in\mn\}$ is a base at $\top$. Therefore, $(L, \sigma (L))$ is not first countable. One can easily check that $(L, \sigma (L))$ is sober.
\end{example}

\section{$\omega$-Rudin-sets and $\omega$-well-filtered determined sets}

In this section, based on the topological Rudin's Lemma, we introduce and study two new classes of closed subsets in $T_0$ spaces - $\omega$-Rudin sets and  $\omega$-well-filtered determined closed sets lying between the class of all closures of countable directed subsets and that of irreducible closed subsets.

For a $T_0$ space $X$ and $\mathcal{K}\subseteq \mathord{\mathsf{K}}(X)$, let $M(\mathcal{K})=\{A\in \mathcal C(X) : K\bigcap A\neq\emptyset \mbox{~for all~} K\in \mathcal{K}\}$ (that is, $\mathcal A\subseteq \Diamond A$) and $m(\mathcal{K})=\{A\in \mathcal C(X) : A \mbox{~is a minimal menber of~} M(\mathcal{K})\}$.

\begin{definition}\label{rudinset}
		Let $X$ be a $T_0$ space. A nonempty subset  $A$  of $X$  is said to have the \emph{$\omega$-Rudin property}, if there exists a countable filtered family $\mathcal K\subseteq \mathord{\mathsf{K}}(X)$ such that $\overline{A}\in m(\mathcal K)$ (that is,  $\overline{A}$ is a minimal closed set that intersects all members of $\mathcal K$). Let $\mathsf{RD}_\omega(X)=\{A\in \mathcal C(X) : A\mbox{~has $\omega$-Rudin property}\}$. The sets in $\mathsf{RD}_\omega(X)$ will also be called \emph{$\omega$-Rudin sets}.
\end{definition}

\begin{lemma}\label{rudinimage}
	Let $X, Y$ be two $T_0$ spaces and $f : X\longrightarrow Y$ a continuous mapping. If $A\in \mathsf{RD}_\omega(X)$, then $\overline{f(A)}\in \mathsf{RD}_\omega(Y)$.
\end{lemma}
\begin{proof} Since $A\in \mathsf{RD}_\omega(X)$, there exists a countable filtered family $\mathcal K\subseteq \mathord{\mathsf{K}}(X)$ such that $A\in m(\mathcal K)$. Let $\mathcal{K}_f=\{\ua f(K\cap A) : K\in \mathcal K\}$. Then $\mathcal{K}_f\subseteq \mathord{\mathsf{K}}(Y)$ is countable filtered. For each $K\in \mathcal K$, since $K\cap A\neq\emptyset$, we have $\emptyset\neq f(K\cap A)\subseteq \ua f(K\cap A)\cap \overline{f(A)}$. So $\overline{f(A)}\in M(\mathcal{K}_f)$. If $B$ is a closed subset of $\overline{f(A)}$ with $B\in M(\mathcal{K}_f)$, then $B\cap\ua f(K\cap A)\neq\emptyset$ for every $K\in \mathcal K$. So $K\cap A\cap f^{-1}(B)\neq\emptyset$ for all $K\in \mathcal K$. It follows that $A=A\cap f^{-1}(B)$ by the minimality of $A$, and consequently, $\overline{f(A)}\subseteq B$. Therefore, $\overline{f(A)}=B$. Thus $\overline{f(A)}\in \mathsf{RD}_\omega(Y)$.
\end{proof}

\begin{proposition}\label{rudinwf}
	Let $X$ be a $T_0$ space and  $Y$ an $\omega$-well-filtered space. If $f : X\longrightarrow Y$ is continuous and $A\in\mathsf{RD}_\omega(X)$, then there exists a unique $y_A\in X$ such that $\overline{f(A)}=\overline{\{y_A\}}$.
\end{proposition}
\begin{proof}
Since $A\in\mathsf{RD}_\omega(X)$, there exists a countable filtered family $\mathcal K\subseteq \mathord{\mathsf{K}}(X)$ such that $\overline{A}\in m(\mathcal K)$. Let $\mathcal{K}_f=\{\ua f(K\cap \overline{A}) : K\in \mathcal K\}$. Then $\mathcal{F}_f\subseteq \mathord{\mathsf{K}}(Y)$ is countable filtered. By the proof of Lemma \ref{rudinimage}, $\overline{f(A)}\in m(\mathcal{K}_f)$. Since $Y$ is $\omega$-well-filtered, we have $\bigcap_{K\in \mathcal{K}}\ua f(K\cap \overline{A})\cap \overline{f(A)}\neq\emptyset$. Select a $y_A\in \bigcap_{K\in \mathcal K} \ua f(K\cap \overline{A})\cap \overline{f(A)}$. Then $\overline{\{y_A\}}\subseteq \overline{f(A)}$ and $K\cap \overline{A}\cap f^{-1}(\overline{\{y_A\}})\neq\emptyset$ for all $K\in \mathcal K$. It follows that $\overline{A}=\overline{A}\cap f^{-1}(\overline{\{y_A\}})$ by the minimality of $\overline{A}$, and consequently, $\overline{f(A)}\subseteq \overline{\{y_A\}}$. Therefore, $\overline{f(A)}=\overline{\{y_A\}}$. The uniqueness of $y_A$ follows from the $T_0$ separation of $Y$.
\end{proof}

Motivated by Proposition \ref{rudinwf}, we introduce the following concept.

\begin{definition}\label{WDspace}
	 A subset $A$ of a $T_0$ space $X$ is called an \emph{$\omega$-well-filtered determined set}, $\wdd_\omega$ \emph{set} for short, if for any continuous mapping $ f:X\longrightarrow Y$
to an $\omega$-well-filtered space $Y$, there exists a unique $y_A\in Y$ such that $\overline{f(A)}=\overline{\{y_A\}}$.
Denote by $\mathsf{WD}_\omega(X)$ the set of all closed $\omega$-well-filtered determined subsets of $X$.
\end{definition}

Obviously, a subset $A$ of a space $X$ is $\omega$-well-filtered determined if{}f $\overline{A}$ is $\omega$-well-filtered determined.

\begin{proposition}\label{DRWIsetrelation}
	Let $X$ be a $T_0$ space. Then $\mathcal{S}_c(X)\subseteq \mathcal{D}_c^\omega(X)\subseteq \mathsf{RD}_\omega(X)\subseteq\mathsf{WD}_\omega(X)\subseteq\ir_c(X)$.
\end{proposition}
\begin{proof} Obviously, $\mathcal{S}_c(X)\subseteq \mathcal{D}_c^\omega(X)$. Now we prove that the closure of a countable directed subset $D$ of $X$ is an $\omega$-Rudin set. Let
$\mathcal K_D=\{\ua d : d\in D\}$. Then $\mathcal K_D\subseteq \mathord{\mathsf{K}}(X)$ is countable filtered and $\overline{D}\in M(\mathcal K_D)$. If $A\in M(\mathcal K_D)$, then $d\in A$ for every $d\in D$, and hence $\overline{D}\subseteq A$. So $\overline{D}\in m(\mathcal K_D)$. Therefore $\overline{D}\in \mathsf{RD}_\omega(X)$. By Proposition \ref{rudinwf}, $\mathsf{RD}_\omega(X)\subseteq\mathsf{WD}_\omega(X)$. Finally we show $\mathsf{WD}_\omega(X)\subseteq\ir_c(X)$. Let $A\in\mathsf{WD}_\omega(X)$.
	Since $\eta_X: X\longrightarrow X^s,\ x\mapsto\da x$, is a continuous mapping to an $\omega$-well-filtered space ($X^s$ is sober), there exists $C\in \ir_c(X)$ such that $\overline{\eta_X(A)}=\overline{\{C\}}$.
	Let $U\in\mathcal O(X)$.
	Note that $$\begin{array}{lll}
	A\cap U\neq\emptyset &\Leftrightarrow& \eta_X(A)\cap  \Diamond U\neq\emptyset\\
	&\Leftrightarrow&\{C\}\cap \Diamond U\neq\emptyset\\
	&\Leftrightarrow& C\in \Diamond U\\
	&\Leftrightarrow& C\cap U\neq\emptyset.
	\end{array}$$
	This implies that $A=C$, and hence $A\in \ir_c(X)$.
\end{proof}

\begin{lemma}\label{WDimage}
Let $X,Y$ be two $T_0$ spaces. If $f:X\longrightarrow Y$ is a continuous mapping and $A\in \wdd_\omega (X)$, then $\overline{f(A)}\in \wdd_\omega (Y)$.
\end{lemma}
\begin{proof}	Let $Z$ is an $\omega$-well-filtered space and $g:Y\longrightarrow Z$ is a continuous mapping.
Since $g\circ f:X\longrightarrow Z$ is continuous and $A\in \wdd_\omega (X)$, there is $z\in Z$ such that $\overline{g(\overline{f(A)})}=\overline{g\circ f(A)}=\overline{\{z\}}$. Thus $\overline{f(A)}\in \wdd_\omega (Y)$.
\end{proof}

\begin{lemma}\label{Rudinsetprod}
	Let	$X=\prod_{i<\omega}X_i$ be the product of a countable family $\{X_i:i<\omega\}$ of $T_0$ spaces and $A\in\ir_c (X)$. Then the following conditions are equivalent:
\begin{enumerate}[\rm (1)]
	\item $A$ is an $\omega$-Rudin set.
	\item $p_i(A)$ is an $\omega$-Rudin set for each $i\in I$.
\end{enumerate}
\end{lemma}
\begin{proof}(1) $\Rightarrow$ (2):  By Corollary \ref{irrcprod} and Lemma \ref{rudinimage}.
	
	(2) $\Rightarrow$ (1): For each $i <\omega$, by Corollary \ref{irrcprod}, $p_i(A)\in \ir_c(X_i)$, and hence by condition (2), there is a countable filtered family $\mathcal K_i\in \mk (X_i)$ such that $p_i(A)\in m(\mathcal K_i)$. Let $\mathcal K=\{K_{\varphi}=\prod_{i<\omega}{\varphi(i)}: \varphi\in\prod_{i<\omega}\mathcal K_i\}$. Then by Tychonoff's Theorem  (see \cite[pp.184, 3.2.4]{Engelking}), $\mathcal K\subseteq \mathord{\mathsf{K}}(\prod_{i<\omega}X_i)$ and $\mathcal K$ is countable filtered because all $\mathcal K_i$ are countable filtered. $\forall \varphi\in\prod_{i<\omega}\mathcal K_i, i<\omega$, ${\varphi(i)}\cap p_i(A)\neq\emptyset$ since $p_i(A)\in m(\mathcal K_i)$. So by Lemma \ref{irrprod} and Corollary \ref{irrcprod} we have
	$$K_{\varphi}\cap A=\left(\prod_{i<\omega}\varphi(i)\right)\bigcap\left(\prod_{i<\omega} p_i(A)\right)=\prod_{i<\omega}\left({\varphi(i)}\cap p_i(A)\right)\neq\emptyset.$$
It follows that $A\in M(\mathcal K)$. Now we show that $A\in m(\mathcal K)$. If $B$ is a closed subset of $A$ that meets all members of $\mathcal K$, then by Lemma \ref{t Rudin}, $B$ contains an minimal irreducible closed subset $C$ that still meets all
members of $\mathcal{K}$. Then for each $i <\omega$, $p_i(C)\in M(\mathcal K_i)$, so $p_i(C)=p_i(A)$ by $p_i(A)\in m(\mathcal K_i)$ and $p_i(C)\subseteq p_i(A)$. By Corollary \ref{irrcprod} again, we have $A=\prod_{i<\omega}p_i(A)=\prod_{i<\omega}p_i(C)=C$, and hence $A=C$. Thus $A\in \mathsf{RD}_\omega(\prod_{i<\omega}X_i)$. 	
\end{proof}

\begin{lemma}\label{WDsetprod}
	Let	$\{X_i: 1\leq i\leq n\}$ be a finite family of $T_0$ spaces and $X=\prod\limits_{i=1}^{n}X_i$ the product space. For $A\in\ir_c (X)$, the following conditions are equivalent:
\begin{enumerate}[\rm (1)]
	\item $A\in \wdd_\omega (X)$.
	\item $p_i(A)\in \wdd_\omega (X_i)$ for each $1\leq i\leq n$.
\end{enumerate}
\end{lemma}

\begin{proof} (1) $\Rightarrow$ (2): By Corollary \ref{irrcprod} and Lemma \ref{WDimage}.

(2) $\Rightarrow$ (1): By induction, we need only to prove the implication for the case of $n=2$. Let $A_1=p_1(A)$ and $A_2= p_2(A)$. Then by condition (2), $(A_1, A_2)\in \mathcal \wdd_\omega (X_1)\times \wdd_\omega (X_2)$ (note that $(A_1, A_2)\in \ir_c (X_1)\times \ir_c (X_2)$ by Corollary \ref{irrcprod}). Now we show that the product $A_1\times A_2\in\wdd_\omega (X)$. Let $f : X_1\times X_2 \longrightarrow Y$ a continuous mapping from $X_1\times X_2$ to an $\omega$-well-filtered space $Y$. For each $b\in X_2$, $X_1$ is homeomorphic to $X_1\times \{b\}$ (as a subspace of $X_1\times X_2$) via the homeomorphism $\mu_b : X_1 \longrightarrow X_1\times \{b\}$ defined by $\mu_b(x)=(x, b)$. Let $i_b : X_1\times \{b\}\longrightarrow X_1\times X_2$ be the embedding of $X_1\times \{b\}$ in $X_1\times X_2$. Then  $f_{b}=f\circ i_b \circ \mu_b : X_1 \longrightarrow Y$, $f_b(x)=f((x, b))$,  is continuous. Since $A_1\in \mathcal \wdd_\omega (X_1)$, there is a unique $y_b\in Y$ such that $\overline{f(A_1\times \{b\})}=\overline{f_b(A_1)}=\overline{\{y_b\}}$. Define a mapping $g_A : X_2 \longrightarrow Y$ by $g_A(b)=y_b$. For each $V\in \mathcal O(Y)$,
$$\begin{array}{lll}
	g_A^{-1}(V)& =\{b\in X_2 : g_A(b)\in V\}\\
	           & =\{b\in X_2 : \overline{f_b(A_1)}\cap V\neq\emptyset\}\\
	           & =\{b\in X_2 : \overline{f(A_1\times \{b\})}\cap V\neq\emptyset\}\\
	           & =\{b\in X_2 : f(A_1\times \{b\})\cap V\neq\emptyset\}\\
               & =\{b\in X_2 : (A_1\times \{b\})\cap f^{-1}(V)\neq\emptyset\}.\\
	\end{array}$$
 Therefore, for each $b\in g_A^{-1}(V)$, there is an $a_1\in A_1$ such that $(a_1, b)\in f^{-1}(V)\in \mathcal O(X_1\times X_2)$, and hence there is $(U_1, U_2)\in \mathcal O(X_1)\times \mathcal O(X_2)$ such that $(a_1, b)\in U_1\times U_2\subseteq  f^{-1}(V)$. It follows that $b\in U_2\subseteq g_A^{-1}(V)$. Thus $g_A : X_2 \longrightarrow Y$ is continuous. Since $A_2\in \mathcal \wdd_\omega (X_2)$, there is a unique $y_A\in Y$ such that $\overline{g_A(A_2)}=\overline{\{y_A\}}$. Therefore, by Lemma \ref{irrprod} and Corollary \ref{irrcprod}, we have
 $$\begin{array}{lll}
     \overline{f(A)} & =\overline{f(A_1\times A_2)}\\
	                       & =\overline{\bigcup\limits_{a_2\in A_2}f(A_1\times \{a_2\})}\\
	                       & =\overline{\bigcup\limits_{a_2\in A_2}\overline{f(A_1\times \{a_2\})}}\\
	                       & =\overline{\bigcup\limits_{a_2\in A_2}\overline{\{g_A(a_2)\}}}\\
                           & =\overline{\bigcup\limits_{a_2\in A_2}\{g_A(a_2)\}}\\
                           & =\overline{g_A(A_2)}\\
                           & =\overline{\{y_A\}}.\\
	\end{array}$$
Thus $A\in \wdd_\omega (X)$.
\end{proof}

By Corollary \ref{irrcprod} and Lemma \ref{WDsetprod}, we get the following result.

\begin{corollary}\label{WDclosedsetprod}
	Let	$X=\prod\limits_{i=1}^{n}X_i$ be the product of a finitely family $\{X_i: 1\leq i\leq n\}$ of $T_0$ spaces. If $A\in\wdd_\omega (X)$, then $A=\prod\limits_{i=1}^{n}p_i(X_i)$ and $p_i(A)\in \wdd_\omega (X_i)$ for all $1\leq i \leq n$.
\end{corollary}

\begin{question}\label{WDinfiniteset-prodquestion} Let $X=\prod_{i<\omega}X_i$ be the product space of a countable family $\{X_i: i<\omega\}$ of $T_0$ spaces. If all $A_i\subseteq X_i$ $(i<\omega )$ are $\omega$-$\wdd$ sets, must  the product set $\prod_{i<\omega}A_i$ be an $\omega$-$\wdd$ set of $X$?
\end{question}

\begin{theorem}\label{soberequiv} For a $T_0$ space $X$, the following conditions are equivalent:
	\begin{enumerate}[\rm (1)]
		\item $X$ is sober.
		\item $X$ is an $\omega$-$d$-space and $\ir_c(X)=\mathsf D_c^\omega (X)$.
        \item $X$ is $\omega$-well-filtered and $\ir_c(X)=\mathsf D_c^\omega (X)$.
		\item $X$ is $\omega$-well-filtered and $\ir_c(X)=\kf_\omega (X)$.
		\item $X$ is $\omega$-well-filtered and $\ir_c(X)=\wdd_\omega (X)$.
	\end{enumerate}
\end{theorem}
\begin{proof}
	By Proposition \ref{d-spacecharac1}, Proposition \ref{wfismc} and Proposition \ref{DRWIsetrelation}, we only need to check (5) $\Rightarrow$ (1).
	Assume $X$ is  $\omega$-well-filtered and $\ir_c(X)=\wdd_\omega (X)$. Let $A\in\ir_c(X)$. Since the identity $id_X : X\longrightarrow X$ is continuous, there is a unique $x\in X$ such that $\overline{A}=\overline{\{x\}}$. So $X$ is sober.
\end{proof}

\begin{example}\label{examp1}
	Let $X$ be a countable infinite set and endow $X$ with the cofinite topology (having the complements of the finite sets as open sets). The
resulting space is denoted by $X_{cof}$. Then $\mk (X_{cof})=2^X\setminus \{\emptyset\}$ (that is, all nonempty subsets of $X$), and hence $X_{cof}$ is a locally compact and first countable $T_1$ space. Let $\mathcal K=\{X\setminus F : F\in X^{(<\omega)}\}$. It is easy to check that $\mathcal K\subseteq \mk (X_{cof})$ is countable filtered and $X\in m(\mathcal K)$. Therefore, $X\in \kf_\omega(X)$ but $X\not\in \md_c(X)$, and hence $X\not\in \md_c^\omega (X)$.  Thus $\kf_\omega(X)\neq \md_c^\omega(X)$ and $\wdd_\omega(X)\neq \md_c^\omega(X)$. $X_{cof}$ is a $d$-space. Since $X\in \ir_c (X_{cof})\setminus \mathcal S_c(X_{cof})$, $X_{cof}$ is non-sober, and hence is not $\omega$-well-filtered by Theorem \ref{maintheorem2}. In fact, $\mathcal K=\{X\setminus F : F\in X^{(<\omega)}\}\subseteq \mk (X_{cof})$ is countable filtered and $\bigcap \mathcal K=X\setminus \bigcup X^{(<\omega)}=X\setminus X=\emptyset$, but $X\setminus F\neq\emptyset$ for all $F\in X^{(<\omega)}$.
\end{example}

\begin{example}\label{examp2}
	Let $L$ be the complete lattice constructed by Isbell \cite{isbell}. Then it is not sober, and by Corollary \ref{xi-lawson2}, $\Sigma L$ is a well-filtered space, and hence an $\omega$-well-filtered. By Theorem \ref{soberequiv}, $\ir_c(X)\neq\kf_\omega(X)$ and $\ir_c(X)\neq\wdd_\omega (X)$.
\end{example}

\begin{question}\label{R-Wquestion} Does $\mathsf{RD}_\omega(X)=\mathsf{WD}_\omega(X)$ hold for ever $T_0$ space $X$?
\end{question}

\section{$\omega$-well-filtered reflections of $T_0$ spaces}

In this section, we present a direct construction of the  $\omega$-well-filtered reflections of $T_0$ spaces. Some basic properties of $\omega$-well-filtered reflections of $T_0$ spaces are investigated.

\begin{definition}\label{WFtion}
	Let $X$ be a $T_0$ space. An \emph{$\omega$-well-filtered reflection} of $X$ is a pair $\langle \widetilde{X}, \mu\rangle$ consisting of an $\omega$-well-filtered space $\widetilde{X}$ and a continuous mapping $\mu :X\longrightarrow \widetilde{X}$ satisfying that for any continuous mapping $f: X\longrightarrow Y$ to an $\omega$-well-filtered space, there exists a unique continuous mapping $f^* : \widetilde{X}\longrightarrow Y$ such that $f^*\circ\mu=f$, that is, the following diagram commutes.\\
\begin{equation*}
	\xymatrix{
		X \ar[dr]_-{f} \ar[r]^-{\mu}
		&\widetilde{X}\ar@{.>}[d]^-{f^*}\\
		&Y}
	\end{equation*}

\end{definition}

$\omega$-well-filtered reflections, if they exist, are unique up to homeomorphism. We shall use $X^{\omega\mbox{-}w}$ to denote the space of the $\omega$-well-filtered reflection of $X$ if it exists.

	Let $X$ be a $T_0$ space. Then by Proposition \ref{DRWIsetrelation}, $\wdd_\omega (X)\subseteq \ir_c(X)$, and whence the space $P_H(\wdd_\omega(X))$ has the topology $\{\Diamond U : U\in \mathcal O(X)\}$, where
$\Diamond U=\{A\in \wdd_\omega(X) : A\cap U\neq\emptyset\}$. The closed subsets of $P_H(\wdd_\omega(X))$ are exactly the set of forms $\Square C=\downarrow_{\wdd_\omega(X)}C$ with $C\in\mathcal C(X)$.

\begin{lemma}\label{lemmaclosure}
	Let $X$ be a $T_0$ space and $A\subseteq X$. Then $\overline{\eta_X(A)}=\overline{\eta_X\left(\overline{A}\right)}=\overline{\Box A}=\Box \overline{A}$ in $P_H(\wdd_\omega(X))$.
\end{lemma}
\begin{proof}
	Clearly, $\eta_X(A)\subseteq \Box A\subseteq \Box\overline{A}$, $\eta_X\left(\overline{A}\right)\subseteq \Box\overline{A}$ and $\Box\overline{A}$ is closed in $P_H(\wdd_\omega(X))$. It follows that
	$$\overline{\eta_X(A)}\subseteq \overline{\Box A}\subseteq \Box \overline{A}\ \text{ and }\ \overline{\eta_X(A)}\subseteq\overline{\eta_X\left(\overline{A}\right)}\subseteq\Box\overline{A}.$$
	To complete the proof, we need to show $\Box\overline{A}\subseteq \overline{\eta_X(A)}$.
	Let $F\in \Box\overline{A}$. Suppose $U\in\mathcal O(X)$ such that $F\in\Diamond U$, that is, $F\cap U\neq\emptyset$. Since $F\subseteq \overline{A}$, we have $A\cap U\neq\emptyset$. Let $a\in A\cap U$. Then $\da a\in \Diamond U\cap \eta_X(A)\neq\emptyset$. This implies that $F\in \overline{\eta_X(A)}$. Whence $\Box\overline{A}\subseteq \overline{\eta_X(A)}$.	
\end{proof}

\begin{lemma}\label{lemmaeta}
	The mapping $\eta_X:X\longrightarrow P_H(\wdd_\omega(X))$ defined by
	$$\forall x\in X, \ \eta_X(x)=\da x,$$
	is a topological embedding.
\end{lemma}
\begin{proof}
For $U\in\mathcal O(X)$, we have $$\eta_X^{-1}(\Diamond U)=\{x\in X: \da x\in\Diamond U\}=\{x\in X: x\in U\}=U,$$ so $\eta_X$ is continuous.
In addition, we have $$
\eta_X(U)=\{\da x: x\in U\}
=\{\da x: \da x\in\Diamond U\}
=\Diamond U\cap \eta_X(X),$$
which implies that $\eta_X$ is an open mapping to $\eta_X(X)$, as a subspace of $P_H(\wdd_\omega(X))$.
As $\eta_X$ is an injection, $\eta_X$ is a topological embedding.
\end{proof}

\begin{lemma}\label{lemmaWDirr}
	Let $X$ be a $T_0$ space and $A$ a nonempty subset of $X$. Then the following conditions are equivalent:
	\begin{enumerate}[\rm (1)]
		\item $A$ is irreducible in $X$.
		\item $\Box A$ is irreducible in $P_H(\wdd_\omega(X))$.
        \item $\Box \overline{A}$ is irreducible in $P_H(\wdd_\omega(X))$.
	\end{enumerate}
\end{lemma}
\begin{proof}
	(1) $\Rightarrow$ (3): Assume $A$ is irreducible. Then $\eta_X(A)$ is irreducible in $P_H(\wdd_\omega(X))$ by Lemma \ref{irrimage} and Lemma \ref{lemmaeta}. By Lemma \ref{irrsubspace} and Lemma \ref{lemmaclosure}, $\Box \overline{A}=\overline{\eta_X(A)}$ is irreducible in $P_H(\wdd_\omega(X))$.
	
	(3) $\Rightarrow$ (1): Assume $\Box \overline{A}$ is irreducible. Let $A\subseteq B\cup C$ with  $B,C\in\mathcal C(X)$. By Proposition \ref{DRWIsetrelation}, $\wdd_\omega (X)\subseteq \ir_c(X)$, and consequently, we have $\Box\overline{A}\subseteq \Box B\cup\Box C$. Since $\Box\overline{A}$ is irreducible,  $\Box\overline{A}\subseteq \Box B$ or $\Box\overline{A}\subseteq C$, showing that  $\overline{A}\subseteq B$ or $\overline{A}\subseteq C$, and consequently, $A\subseteq B$ or $A\subseteq C$, proving $A$ is irreducible.

    (2) $\Leftrightarrow$ (3): By Lemma \ref{irrsubspace} and Lemma \ref{lemmaclosure}.

\end{proof}

\begin{lemma}\label{lemmafstar}
Let $X$ be a $T_0$ space and $f:X\longrightarrow Y$ a continuous mapping from $X$ to an $\omega$-well-filtered space $Y$. Then there exists a unique continuous mapping $f^* :P_H(\wdd_\omega(X))\longrightarrow Y$ such that $f^*\circ\eta_X=f$, that is, the following diagram commutes.
\begin{equation*}
\xymatrix{
	X \ar[dr]_-{f} \ar[r]^-{\eta_X}
	&P_H(\wdd_\omega(X))\ar@{.>}[d]^-{f^*}\\
	&Y}
\end{equation*}	
\end{lemma}
\begin{proof}For each $A\in\wdd_\omega(X)$, there exists a unique $y_A\in Y$ such that $\overline{f(A)}=\overline{\{y_A\}}$. Then we can define a mapping $f^*:P_H(\wdd_\omega(X))\longrightarrow Y$ by
$$\forall A\in\wdd_\omega(X),\ \ f^*(A)=y_A.$$

{Claim 1:}  $f^*\circ \eta_X=f$.

Let $x\in X$. Since $f$ is continuous, we have
$\overline{f\left(\overline{\{x\}}\right)}=\overline{f(\{x\})}=\overline{\{f(x)\}}$,
so
$f^*\left(\overline{\{x\}}\right)=f(x)$. Thus $f^*\circ \eta_X=f$.

{Claim 2:}  $f^*$ is continuous.

Let $V\in\mathcal O(Y)$. Then
$$\begin{array}{lll}
(f^*)^{-1}(V)&=&\{A\in\wdd_\omega(X): f^*(A)\in V\}\\
&=&\{A\in\wdd_\omega(X): \overline{\{f^*(A)\}}\cap V\neq\emptyset\}\\
&=&\{A\in\wdd_\omega(X): \overline{f(A)}\cap V\neq\emptyset\}\\
&=&\{A\in\wdd_\omega(X): f(A)\cap V\neq\emptyset\}\\
&=&\{A\in\wdd_\omega(X): A\cap f^{-1}(V)\neq\emptyset\}\\
&=&\Diamond f^{-1}(V),
\end{array}$$
which shows that $(f^*)^{-1}(V)$ is open in $P_H(\wdd_\omega(X))$. Thus  $f^*$ is continuous.

{Claim 3:}  The mapping $f^*$ is unique such that $f^*\circ \eta_X=f$.

Assume $g:P_H(\wdd_\omega(X))\longrightarrow Y$ is a continuous mapping such that $g\circ\eta_X=f$.  Let $A\in\wdd_\omega(X)$. We need to show $g(A)=f^*(A)$.
Let $a\in A$.  Then $\overline{\{a\}}\subseteq A$, implying that $g(\overline{\{a\}})\leq_Y g(A)$, that is,  $g\left(\overline{\{a\}}\right)=f(a)\in\overline{\{g(A)\}}$. Thus $\overline{\{f^*(A)\}}=\overline{f(A)}\subseteq \overline{\{g(A)\}}$.
In addition, since $A\in\overline{\eta_X(A)}$ and $g$ is continuous, $g(A)\in g\left(\overline{\eta_X(A)}\right)\subseteq\overline{g(\eta_X(A))}=\overline{f(A)}=\overline{\{f^*(A)\}}$, which implies that $\overline{\{g(A)\}}\subseteq \overline{\{f^*(A)\}}$.  So $\overline{\{g(A)\}}=\overline{\{f^*(A)\}}$. Since $Y$ is $T_0$, $g(A)=f^*(A)$. Thus $g=f^*$.
\end{proof}

\begin{lemma}\label{lemmaBoxwd}
	Let $X$ be a $T_0$ space and $C\in\mathcal C(X)$. Then the following conditions are equivalent:
	\begin{enumerate}[\rm (1)]
		\item $C$ is $\omega$-well-filtered determined in $X$.
		\item $\Box C$ is $\omega$-well-filtered determined in $P_H(\wdd_\omega(X))$.
	\end{enumerate}
\end{lemma}
\begin{proof}
(1) $\Rightarrow$ (2): By Propositions \ref{WDimage}, Lemma \ref{lemmaclosure} and Lemma \ref{lemmaeta}.

(2) $\Rightarrow$ (1).  Let $Y$ be an $\omega$-well-filtered space and $f:X\longrightarrow Y$  a continuous mapping. By Lemma \ref{lemmafstar}, there exists a continuous mapping $f^* :P_H(\wdd_\omega(X))\longrightarrow Y$ such that $f^*\circ\eta_X=f$.
Since $\Box C=\overline{\eta_X(C)}$ is $\omega$-well-filtered determined and $f^*$ is continuous,  there exists a unique $y_C\in Y$ such that
$\overline{f^*\left(\overline{\eta_X(C)}\right)}=\overline{\{y_C\}}$. Furthermore, we have
$$\overline{\{y_C\}}=\overline{f^*\left(\overline{\eta_X(C)}\right)}=\overline{f^*(\eta_X(C))}=\overline{f(C)}.$$
So $C$ is $\omega$-well-filtered determined.
\end{proof}

\begin{theorem}\label{WDwf}
	Let $X$ be a $T_0$ space. Then $P_H(\wdd_\omega(X))$ is an $\omega$-well-filtered space.
\end{theorem}
\begin{proof}
	Since $X$ is $T_0$, one can deduce that $P_H(\wdd_\omega(X))$ is $T_0$. Let $\{
	\mathcal K_i:i\in I\}\subseteq\mk(P_H(\wdd_\omega(X)))$ be a countable filtered family and $U\in\mathcal O(X)$ such that
	$\bigcap_{i\in I}\mathcal K_i\subseteq \Diamond U$. We need to show $\mathcal K_i\subseteq\Diamond U$ for some $i\in I$.
	Assume, on the contrary, $\mathcal K_i\nsubseteq \Diamond U$, i.e., $\mathcal K_i\cap\Box(X\setminus U)\neq\emptyset$,  for any $i\in I$.
	
	Let $\mathcal A=\{C\in\mathcal C(X): C\subseteq X\setminus U \text{ and } \mathcal K_i\cap \Box C\neq\emptyset \text{ for all }  i\in I\}$. Then we have the following two facts.
	
	{(a1)} $\mathcal A\neq\emptyset$ because $X\setminus U\in\mathcal A$.
	
	{(a2)} For any filtered family $\mathcal F\subseteq\mathcal A$, $\bigcap\mathcal F\in\mathcal A$.
	
	Let $F=\bigcap\mathcal F$. Then $F\in \mathcal C(X)$ and $F\subseteq X\setminus U$. Assume, on the contrary, $F\notin\mathcal A$. Then there exists $i_0\in I$ such that $\mathcal K_{i_0}\cap \Box F=\emptyset$. Note that $\Box F=\bigcap_{C\in\mathcal F}\Box C$, implying that $\mathcal K_{i_0}\subseteq\bigcup_{C\in\mathcal F}\Diamond (X\setminus C)$ and $\{\Diamond (X\setminus C) : C\in\mathcal F\}$ is a directed family since $\mathcal F$ is filtered. Then there is $C_0\in\mathcal F$ such that $\mathcal K_{i_0}\subseteq \Diamond (X\setminus C_0)$, i.e., $\mathcal K_{I_0}\cap\Box C_0=\emptyset$,  contradicting $C_0\in\mathcal A$. Hence $F\in\mathcal A$.
	
	By Zorn's Lemma, there exists a minimal element $C_m$ in $\mathcal A$ such that $\Box C_m$ intersects all members of $\mathcal K$. Clearly,  $\Box C_m$ is also a minimal closure set that  intersects all members of $\mathcal K$, hence is an $\omega$-Rudin set in $P_H(\wdd_\omega(X))$. By Proposition \ref{DRWIsetrelation} and Lemma \ref{lemmaBoxwd}, $C_m$ is $\omega$-well-filtered determined. So $C_m\in \Box C_m\cap\bigcap\mathcal K\neq\emptyset$. It follows that $\bigcap\mathcal K\nsubseteq \Diamond(X\setminus C_m)\supseteq \Diamond U$, which implies that $\bigcap\mathcal K\nsubseteq \Diamond U$, a contradiction.
\end{proof}

By Lemma \ref{lemmafstar} and Theorem \ref{WDwf}, we have the following result.

\begin{theorem}\label{WFilterification}
	Let $X$ be a $T_0$ space and $X^{\omega\mbox{-}w}=P_H(\wdd_\omega(X))$. Then the pair $\langle X^{\omega\mbox{-}w}, \eta_X\rangle$, where $\eta_X :X\longrightarrow X^{\omega\mbox{-}w}$, $x\mapsto\overline{\{x\}}$, is the $\omega$-well-filtered reflection of $X$.
\end{theorem}

\begin{corollary}\label{WFreflective}
	The category $\mathbf{Top}_{\omega\mbox{-}w}$ of all $\omega$-well-filtered spaces is a reflective full subcategory of  $\mathbf{Top}_0$.
\end{corollary}

\begin{corollary}\label{WFfuctor}
	Let $X,Y$ be two $T_0$ spaces and $f:X\longrightarrow Y$  a continuous mapping. Then there exists a unique continuous mapping $f^{\omega\mbox{-}w}:X^{\omega\mbox{-}w}\longrightarrow Y^{\omega\mbox{-}w}$ such that $f^{\omega\mbox{-}w}\circ \eta_X=\eta_Y\circ f$, that is, the following diagram commutes.
		\begin{equation*}
	\xymatrix{
		X \ar[d]_-{f} \ar[r]^-{\eta_X} &X^{\omega\mbox{-}w}\ar[d]^-{f^{\omega\mbox{-}w}}\\
		Y \ar[r]^-{\eta_Y} &Y^{\omega\mbox{-}w}
	}
	\end{equation*}
For each $A\in \wdd_\omega (X)$, $f^{\omega\mbox{-}w}(A)=\overline{f(A)}$.
\end{corollary}

Corollary \ref{WFfuctor} defines a functor $W_\omega : \mathbf{Top}_0 \longrightarrow \mathbf{Top}_{\omega\mbox{-}w}$, which is the left adjoint to the inclusion functor $I : \mathbf{Top}_{\omega\mbox{-}w} \longrightarrow \mathbf{Top}_0$.

\begin{corollary}\label{WFwdc}
	For a $T_0$ space $X$, the following conditions are equivalent:
	\begin{enumerate}[\rm (1)]
		\item $X$ is $\omega$-well-filtered.
		\item $\mathsf{RD}_\omega(X)=\mathcal S_c(X)$.
        \item $\wdd_\omega (X)=\mathcal S_c(X)$, that is, for each $A\in\wdd_\omega(X)$, there exists a unique $x\in X$ such that $A=\overline{\{x\}}$.
        \item $X\cong X^{\omega\mbox{-}w}$.

	\end{enumerate}
\end{corollary}
\begin{proof}  (1) $\Rightarrow$ (2): Applying Proposition \ref{rudinwf} to the identity $id_X : X \longrightarrow X$.

(2) $\Rightarrow$ (3): By Proposition \ref{DRWIsetrelation}.

(3) $\Rightarrow$ (4): By assumption, $\wdd_\omega(X)=\left\{\overline{\{x\}}:x\in X\right\}$, so $X^{\omega\mbox{-}w}=P_H(\wdd_\omega(X))=P_H(\{\overline {\{x\}} : x\in X\})$, and whence $X\cong X^{\omega\mbox{-}w}$.

(4) $\Rightarrow$ (1): By Theorem \ref{WDwf} or by Proposition \ref{DRWIsetrelation} and Corollary \ref{WFwdc}.
\end{proof}

\begin{remark} By Proposition \ref{d-spacecharac1}, Proposition \ref{DRWIsetrelation} and Corollary \ref{WFwdc}, we can get Proposition \ref{wfismc}.
\end{remark}

\begin{corollary}\label{wfretract}  A retract of an $\omega$-well-filtered space is $\omega$-well-filtered.
\end{corollary}
\begin{proof} Suppose that $Y$ is a retract of an $\omega$-well-filtered space $X$. Then there are continuous mappings $f : X\longrightarrow Y$ and $g : Y\longrightarrow X$ such that $f\circ g=id_Y$. Let $B\in \mathsf{WD}_\omega(Y)$, then by Lemma \ref{WDimage} and Corollary \ref{WFwdc}, there exists a unique $x_B\in X$ such that $\overline{g(B)}=\overline{\{x_B\}}$. Therefore, $B=\overline{f\circ g(B)}=\overline{f(\overline{g(B)})}=\overline{f(\overline{\{x_B\}})}=\overline{\{f(x_B)\}}$. By Corollary \ref{WFwdc}, $Y$ is $\omega$-well-filtered.
\end{proof}

\begin{theorem}\label{wfreflectionprod}
	Let $\{X_i: 1\leq i\leq n\}$ be a finitely family of $T_0$ spaces. Then $(\prod\limits_{i=1}^{n} X_i)^{\omega\mbox{-}w}=\prod\limits_{i=1}^{n}X_{i}^{\omega\mbox{-}w}$ (up to homeomorphism).
\end{theorem}

\begin{proof}	
	Let $X=\prod\limits_{i=1}^{n}X_i$. By Corollary \ref{WDclosedsetprod}, we can define a mapping $\gamma : P_H(\wdd_\omega (X))  \longrightarrow \prod\limits_{i=1}^{n}P_H(\wdd_\omega (X_i))$ by

\begin{center}
$\forall A\in \wdd_\omega (X)$, $\gamma (A)=(p_1(A), p_2(A), ..., p_n(A))$.
\end{center}

By Lemma \ref{WDsetprod} and Corollary \ref{WDclosedsetprod}, $\gamma$ is bijective. Now we show that $\gamma$ is a homeomorphism. For any $(U_1, U_2, ..., U_n)\in \mathcal O(X_1)\times \mathcal O(X_2)\times ... \times \mathcal O(X_n)$, by Lemma \ref{WDsetprod} and Corollary \ref{WDclosedsetprod}, we have

$$\begin{array}{lll}
\gamma^{-1}(\Diamond U_1\times \Diamond U_2\times ... \times\Diamond U_n)&=&\{A\in\wdd_\omega(X): \gamma(A)\in \Diamond U_1\times \Diamond U_2\times ... \times\Diamond U_n\}\\
&=&\{A\in\wdd_\omega(X): p_1(A)\cap U_1\neq\emptyset, p_2(A)\cap U_2\neq\emptyset, ..., p_n(A)\cap U_n\neq\emptyset\}\\
&=&\{A\in\wdd_\omega(X): A\cap U_1\times U_2\times ... \times U_n\neq\emptyset\}\\
&=&\Diamond (U_1\times U_2\times ... \times U_n)\in \mathcal O(P_H(\wdd_\omega (X)), \mbox{~and~}
\end{array}$$

$$\begin{array}{lll}
\gamma (\Diamond (U_1\times U_2\times ... \times U_n))&=&\{\gamma (A): A\in \wdd_\omega (X) \mbox{~and~} A\cap U_1\times U_2\times ... \times U_n\neq\emptyset \}\\
&=&\{\gamma (A): A\in \wdd_\omega (X), \mbox{~and~} p_1(A)\cap U_1\neq\emptyset, ..., p_n(A)\cap U_n\neq\emptyset\}\\
&=&\Diamond U_1\times \Diamond U_2\times ... \times\Diamond U_n\in O(\prod\limits_{i=1}^{n}P_H(\wdd_\omega (X_i))).
\end{array}$$

Therefore, $\gamma : P_H(\wdd_\omega (X))  \longrightarrow \prod\limits_{i=1}^{n}P_H(\wdd_\omega (X_i))$ is a homeomorphism. Thus $X^{\omega\mbox{-}w}$ ($=P_H(\wdd_\omega (X)$) and $\prod\limits_{i=1}^{n}X_i^{\omega\mbox{-}w}$ ($=\prod\limits_{i=1}^{n}P_H(\wdd_\omega (X_i))$ are homeomorphic.
\end{proof}

\begin{question}\label{infiniteWFication} Does $(\prod\limits_{i<\omega }X_i)^w=\prod\limits_{i<\omega}X_i^w$ (up to homeomorphism) hold for any countable family $\{X_i : <\omega\}$ of $T_0$ spaces?
\end{question}

Using $\wdd_\omega$ sets and Corollary \ref{WFwdc}, we show that products of $\omega$-well-filtered spaces are $\omega$-well-filtered.

\begin{theorem}\label{WFprod}
	Let $\{X_i:i\in I\}$ be a family of $T_0$ spaces. Then the following two conditions are equivalent:
	\begin{enumerate}[\rm(1)]
		\item The product space $\prod_{i\in I}X_i$ is $\omega$-well-filtered.
		\item For each $i \in I$, $X_i$ is $\omega$-well-filtered.
	\end{enumerate}
\end{theorem}
\begin{proof}	
	(1) $\Rightarrow$ (2):  For each $i \in I$, $X_i$ is a retract of $\prod_{i\in I}X_i$. By Corollary \ref{wfretract}, $X_i$ is $\omega$-well-filtered.
	
	(2) $\Rightarrow$ (1): Let $X=\prod_{i\in I}X_i$. Suppose $A\in \wdd_\omega (X)$. Then by Corollary \ref{irrcprod}, Proposition \ref{DRWIsetrelation} and Lemma \ref{WDimage}, $A\in \ir_c(X)$ and for each $i \in I$, $p_i(A)\in \wdd_\omega (X_i)$, and consequently, there is a $u_i\in X_i$ such that $p_i(A)=cl_{X_i}\{u_i\}$ by condition (2) and Corollary \ref{WFwdc}. Let $u=(u_i)_{i\in I}$. Then by Lemma \ref{irrprod}, Corollary \ref{irrcprod} and \cite[Proposition 2.3.3]{Engelking}), we have $A=\prod_{i\in I}p_i(A)=\prod_{i\in I}\cl_{u_i}\{u_i\}=\cl_X \{u\}$. Therfore, $X$ is $\omega$-well-filtered by Corollary \ref{WFwdc}.
\end{proof}

\begin{theorem}\label{Wfication=sober}
	For a $T_0$ space $X$, the following conditions are equivalent:
	\begin{enumerate}[\rm (1)]
		\item $X^{\omega\mbox{-}w}$ is the sobrification of $X$, in other words, the $\omega$-well-filtered reflection of $X$ and sobrification of $X$ are the same.
        \item $X^{\omega\mbox{-}w}$ is sober.
		\item $\ir_c(X)=\wdd_\omega(X)$.
	\end{enumerate}
\end{theorem}

\begin{proof}
(1) $\Rightarrow$ (2): Trivial.

(2) $\Rightarrow$ (3): Let $\eta_X^{\omega\mbox{-}w} : X\longrightarrow X^{\omega\mbox{-}w}$ be the canonical topological embedding defined by $\eta_X^{\omega\mbox{-}w}(x)=\overline{\{x\}}$ (see Theorem \ref{WFilterification}). Since the pair $\langle X^s, \eta_{X}^s\rangle$, where $\eta_{X}^s :X\longrightarrow X^s=P_H(\ir_c(X))$, $x\mapsto\overline{\{x\}}$, is the soberification of $X$ and $X^{\omega\mbox{-}w}$ is sober, there exists a unique continuous mapping $(\eta_{X}^{\omega\mbox{-}w})^* :X^s\longrightarrow X^{\omega\mbox{-}w}$ such that $(\eta_{X}^{\omega\mbox{-}w})^*\circ\eta_X^s=\eta_X^{\omega\mbox{-}w}$, that is, the following diagram commutes.
\begin{equation*}
\xymatrix{
	X \ar[dr]_-{\eta_X^{\omega\mbox{-}w}} \ar[r]^-{\eta_X^s}
	&X^s\ar@{.>}[d]^-{(\eta_{X}^{\omega\mbox{-}w})^*}\\
	&X^{\omega\mbox{-}w}}
\end{equation*}	
So for each $A\in\ir_c(X)$, there exists a unique $B\in \wdd_\omega (X)$ such that
$$\downarrow_{\wdd_\omega (X)} A=\overline{\eta_X^{\omega\mbox{-}w}(A)}=\overline{\{B\}}=\downarrow_{\wdd_\omega (X)}B.$$
Clearly, we have $B\subseteq A$. On the other hand, for each $a\in A, \overline {\{a\}}\in \downarrow_{\wdd_\omega (X)} A=\downarrow_{\wdd_\omega (X)}B$, and whence $\overline {\{a\}}\subseteq B$. Thus $A\subseteq B$, and consequently, $A=B$. Thus $A\in \wdd_\omega(X)$.

(3) $\Rightarrow$ (1): If $\wdd_\omega(X)=\ir_c(X)$, then $X^{\omega\mbox{-}w}=P_H(\wdd_\omega (X))=P_H(\ir_c(X))=X^s$, with $\eta_X^{\omega\mbox{-}w}=\eta_X^s : X\longrightarrow X^{\omega\mbox{-}w}$, is the sobrification of $X$.

\end{proof}

\begin{proposition}\label{WFicationcomp}
	A $T_0$ space $X$ is compact if{}f $X^{\omega\mbox{-}w}$ is compact.
\end{proposition}
\begin{proof} By Proposition \ref{DRWIsetrelation}, we have $\mathcal S_c(X)\subseteq \wdd_\omega (X)\subseteq \ir_c (X)$. Suppose that $X$ is compact. For $\{U_i : i\in I\}\subseteq \mathcal O(X)$, if $\wdd_\omega (X)\subseteq \bigcup_{i\in I}\Diamond U_i$, then $X\subseteq \bigcup_{i\in I} U_i$ since $\mathcal S_c(X)\subseteq \wdd_\omega (X)$, and consequently, $X\subseteq \bigcup_{i\in I_0} U_i$ for some $I_0\in I^{(<\omega)}$. It follows that $\wdd_\omega (X)\subseteq \bigcup_{i\in I_0} \Diamond U_i$. Thus $X^{\omega\mbox{-}w}$ is compact. Conversely, if $X^{\omega\mbox{-}w}$ is compact and $\{V_j : j\in J\}$ is a open cover of $X$, then $\wdd_\omega (X)\subseteq \bigcup_{j\in J}\Diamond V_j$. By the compactness of $X^{\omega\mbox{-}w}$, there is a finite subset $J_0\subseteq J$ such that $\wdd_\omega (X)\subseteq \bigcup_{j\in J_0}\Diamond V_j$, and whence $X\subseteq \bigcup_{j\in J_0}V_j$, proving the compactness of $X$.
\end{proof}

Since $\mathcal{S}_c(X)\subseteq \mathsf{WD}_\omega(X)\subseteq\ir_c(X)$ (see Proposition \ref{DRWIsetrelation}), the correspondence $U \leftrightarrow \Diamond_{\wdd_\omega (X)} U$ is a lattice isomorphism between $\mathcal O(X)$ and $\mathcal O(X^{\omega\mbox{-}w})$, and whence we have the following proposition.

\begin{proposition}\label{wficationLHC}
	Let $X$ be a $T_0$ space. Then
	\begin{enumerate}[\rm (1)]
		\item $X$ is locally hypercompact if{}f $X^{\omega\mbox{-}w}$ is locally hypercompact.
        \item $X$ is a C-space if{}f $X^w$ is a C-space.
	\end{enumerate}
\end{proposition}

\begin{proposition}\label{wdficationLC} For a $T_0$ space $X$, the following conditions are equivalent:
\begin{enumerate}[\rm (1)]
		\item $X$ is core compact.
        \item $X^{\omega\mbox{-}w}$ is core compact.
        \item $X^{\omega\mbox{-}w}$ is locally compact.
\end{enumerate}
\end{proposition}
\begin{proof} (1) $\Leftrightarrow$ (2): Since $\mathcal O(X)$ and $\mathcal O(X^{\omega\mbox{-}w})$ are lattice-isomorphic.

(2) $\Rightarrow$ (3): By Theorem \ref{WDwf}, $X^{\omega\mbox{-}w}$ is $\omega$-well-filtered. If $X^{\omega\mbox{-}w}$ is core compact, then $X^{\omega\mbox{-}w}$ is locally compact by Theorem \ref{omegaWFcorecompt-locCompt}.

(3) $\Rightarrow$ (2): Trivial.
\end{proof}

\begin{remark}\label{corecompnotLC} In \cite{Hofmann-Lawson} (see also \cite[Exercise V-5.25]{redbook}) Hofmann and Lawson constructed a core compact $T_0$ space $X$ that is  not locally compact. By Proposition \ref{wdficationLC}, $X^{\omega\mbox{-}w}$ is locally compact. So the local compactness of $X^{\omega\mbox{-}w}$ does not imply the local compactness of $X$.
\end{remark}





\end{document}